\documentclass[12pt]{article}
\usepackage{amsmath}
\usepackage{amsthm}
\usepackage{amssymb}
\usepackage{eucal}
\numberwithin{equation}{section}
\usepackage{graphicx,graphics}
\graphicspath{{figures/}}
\usepackage{tikz}
\usepackage{pgfplots} 
\pgfplotsset{compat=newest}
\usetikzlibrary[pgfplots.groupplots]
\usepgfplotslibrary{groupplots}
\usetikzlibrary{matrix,arrows,decorations.pathmorphing}
\definecolor{k4}{rgb}{0.8,0.8,0.8}
\definecolor{k3}{rgb}{0.6,0.6,0.6}
\definecolor{k2}{rgb}{0.4,0.4,0.4}
\definecolor{k1}{rgb}{0.2,0.2,0.2}

\usepackage{caption}

\usepackage{varwidth}
\DeclareCaptionFormat{myformat}{%
  \begin{varwidth}{\linewidth}%
    \centering
    #1#2#3%
  \end{varwidth}%
}

\numberwithin{equation}{section}

\usepackage{pdfsync}
\usepackage{subcaption}
\def \dis {\displaystyle}

\def \into {\int_\Omega}
\def \confai {-\kern -.5em\rightharpoonup}
\def \cqfd {\hfill$\Box$}

\def\curl{{\rm curl\,}}

\def \ep {\varepsilon}
\def \om {\omega}
\def \Om {\Omega}

\def \NN {\mathbb N}
\def \ZZ {\mathbb Z}
\def \RR {\mathbb R}
\def \CC {\mathbb C}

\def \beq {\begin{equation}}
\def \eeq {\end{equation}}
\def \ba {\begin{array}}
\def \ea {\end{array}}

\def \ecart {\noalign{\medskip}}
\newtheorem{Thm}{Theorem}[section]

\newtheorem{Pro}[Thm]{Proposition}

\newtheorem{Adef}[Thm]{Definition}
\newenvironment{Def}{\begin{Adef}}{\end{Adef}}
\newtheorem{Arem}[Thm]{Remark}
\newenvironment{Rem}{\begin{Arem}}{\end{Arem}}
\newtheorem{Aexa}[Thm]{Example}

\newtheorem{Anot}[Thm]{Notation}

\def \refe #1.{(\ref{#1})}
\def \reff #1.{figure~\ref{#1}}
\def \refs #1.{Section~\ref{#1}}
\def \refss #1.{Subsection~\ref{#1}}
\def \refD #1.{Definition~\ref{#1}}
\def \refT #1.{Theorem~\ref{#1}}
\def \refL #1.{Lemma~\ref{#1}}
\def \refC #1.{Corollary~\ref{#1}}
\def \refP #1.{Proposition~\ref{#1}}
\def \refPt #1.{Properties~\ref{#1}}
\def \refR #1.{Remark~\ref{#1}}
\def \refE #1.{Example~\ref{#1}}
\def \refN #1.{Notation~\ref{#1}}

\title{A strong approximation in $L^2$ for the solutions of the Maxwell system with highly oscillating periodic coefficients. }

\begin{document}
\maketitle
\centerline{\large Juan CASADO-D\'IAZ $^\dag$,\hskip 0.3cm Nourelhouda KHEDHIRI $^{\dag\dag}$}\par\medskip
\centerline{\large Mohamed Lazhar TAYEB $^{\dag\dag}$}
\bigskip \bigskip
\centerline{{$^\dag$} Dpto. de Ecuaciones  Diferenciales y An\'alisis Num\'erico,}
\centerline{Facultad de Matem\'aticas, C. Tarfia s/n}
\centerline{41012 Sevilla, SPAIN}
\bigskip
\centerline{$^{\dag\dag}$ Labo EDP, Department of Mathematics,}
\centerline{Faculty of Sciences of Tunis, El-Manar University,}
\centerline{2092, TUNISIA}
\bigskip
\centerline{e-mail: jcasadod@us.es, khedhirinoorelhouda@gmail.com, mohamedlazhar.tayeb@fst.utm.tn}
\begin{abstract}
We consider a Maxwell system on $\RR^3$ with periodic and highly oscillating coefficients. It is known that the solutions converge in the weak-$\ast$ topology of $L^\infty(0,T;\,L^2(\RR^3))$ to the solution of a similar problem with constant coefficients given as the $H$-limits of the electric permittivity and the magnetic permeability respectively, i.e. the limit in the sense of the homogenization of linear elliptic equations with varying coefficients. However, it is not true that the elliptic corrector also provides a corrector for the solution of the Maxwell system, i.e. an approximation of the solutions in the strong topology of $L^2$. We shall prove that the oscillations in the space variable also produce oscillations in the time variable. We get a corrector consisting of adding to the elliptic corrector the sum of infinitely plane waves in the fast variable. Note that related results have been previously proved for the wave equation. Our proof is based on the two-scale convergence theory for almost periodic functions. One of the novelties is to show how this contains the classical Block decomposition.

\end{abstract}

\textbf{Keywords: } Maxwell system, almost periodic homogenization, corrector result, two-scale convergence.

\medskip

\textbf{2020 Mathematics Subject Classification: }35B27, 78A48.

\section{Introduction.}
The homogenization of the wave equation with varying coefficients is a classical problem, carried out in \cite{CoSp}. Namely, let $\Om$ be a bounded open set of $\RR^N$, $\rho_\ep$ a sequence of functions on $\Om$ uniformly positive and bounded, and $A_\ep$ a sequence of measurable symmetric matrix functions uniformly elliptic and bounded. It is known that, for a subsequence, there exist $\rho$ and $A$ with the same properties, such that for every sequence of source terms
$f_\ep$ which converge weakly in $L^2(0,T;\,L^2(\Om))$ to some $f$ and every sequence of initial data $u_\ep^0$ and $v_\ep^0$ which converge weakly in $H^1_0(\Om)$ and $L^2(\Om)$ respectively to some $u^0$ and $v^0$, the solution $u_\ep$ of (other boundary conditions can be considered)
\beq\label{int1}\left\{\ba{l}\dis \rho_\ep(x)\partial^2_{tt} u_\ep-{\rm div}(A_\ep(x)\nabla u_\ep)=f_\ep\ \hbox{ in }(0,T)\times\Om\\ \ecart \dis u_\ep=0\ \hbox{ on }(0,T)\times\partial\Om\\ \ecart\dis (u_\ep)_{|t=0}=u^0_\ep,\	 \ (\rho_\ep\partial_tu_\ep)_{|t=0}=v_\ep^0\quad\hbox{ in }\Om,\ea\right.\eeq
converge weakly-$\ast$ in $W^{1,\infty}(0,T;H^1_0(\Om))\cap L^\infty(0,T;L^2(\Om))$ to the solution of the analogous problem replacing $\rho_\ep$, $A_\ep$, $f_\ep$, $u^0_\ep$ and $v^0_\ep$ by $\rho$, $A$, $f$, $u^0$ and $v^0$. The function $\rho$ is just the weak-$\ast$ limit of  $\rho_\ep$ in $L^\infty (\Om)$, while $A$ is the $H$-limit of $A_\ep$, i.e. (see e.g. \cite{Mur}, \cite{Spa}, \cite{Tar}) the same matrix that appears in the homogenization of the elliptic equation
\beq\label{int2}-{\rm div}(A_\ep\nabla u_\ep)=f_\ep\ \hbox{ in }\Om,\eeq
where $f_\ep$ is now a sequence which converges strongly in $H^{-1}(\Om)$.
Indeed, (see \cite{CCMM1}) integrating in (\ref{int1}) with respect to $t$ we have that for every $0<t_1<t_2<T$ the sequence
$$\tilde u_\ep(x)=\int_{t_1}^{t_2} u_\ep(t,x)\,dt$$
satisfies a problem similar to (\ref{int2}).\par 
On the other hand, for the homogenization of (\ref{int2}) it is known (\cite{Mur},  \cite{Tar}) that there exists a sequence $H_\ep$ of matrix functions with coefficients in $L^2$ such that denoting by $u$ the weak limit of the solutions of (\ref{int2}) we have that (the convergence holds in $L^2(\Om)^N$ if $u\in W^{1,\infty}(\Om)$)
$$\nabla u_\ep-H_\ep \nabla u\to 0\ \hbox{ in }L^1(\Om)^N.$$
The sequence $H_\ep\nabla u$ is what is known in homogenization as a corrector. Contrarly to $\nabla u$ which only gives an approximation of $\nabla u_\ep$ in a weak topology, $H_\ep \nabla u$ provides an approximation in a strong topology. A similar result also holds for the heat equation.\par

We can ask ourselves if this corrector result also holds for the wave equation, i.e. if the solution $u_\ep$ of (\ref{int1}) satisfies that $\nabla u_\ep-H_\ep\nabla u$ converges strongly to zero at least in $L^1(0,T;\,L^1(\Om))^N.$ Taking a right-hand side independent of $\ep$, it has been proved in \cite{BFM} and \cite{FrMu} that this only holds if the initial condition is \lq\lq well posed\rq\rq which means (in the present case) that ${\rm div}(A_\ep\nabla u_\ep^0)$ and $v_\ep^0/\rho_\ep$ are compact in $H^{-1}(\Om)$ and $L^2(\Om)$ respectively. In general, to get a corrector for $\nabla u_\ep$, a new term must be added to $H_\ep\nabla u$ and it depends non-locally on the initial condition. It will also depend on the right-hand side $f_\ep$ if it only converges weakly in $L^2(0,T;\, L^2(\Om))$. This is due to the fact that the oscillations in the space variable of $A_\ep$ do not only provide oscillations in the space variable for $u_\ep$ but also in the time variables. \par

Assuming that $\rho_\ep=\rho(x/\ep)$, $A_\ep=A(x/\ep)$ with $\rho$, $A$ periodic, and that $\Om=\RR^N$ to avoid the influence of the boundary conditions, a corrector for (\ref{int1}) has been obtained in \cite{BrLe} and \cite{CCMM2}. In \cite{BrLe} the result uses a Block decomposition (see e.g. \cite{BLP}) and considers initial conditions which admit such Block decomposition. The results in \cite{CCMM2} use the two-scale convergence theory for almost periodic functions (\cite{CaGa}, \cite{CaGa2}, \cite{Ngu2}). It assumes initial conditions less general, which are periodic in the fast variable $y=x/\ep$. In \cite{CCMM2}, the right-hand sides can also oscillate and the coefficients are more general including oscillations in the time variable. As a consequence of having a non-local corrector, it is proved   that adding an oscillating damping term to the equation, the limit equation is non-local in general (see 
\cite{BrCa} for related results). \par

We refer to \cite{CCMM3}, \cite{NgLeBr} for the influence of the boundary conditions in the one-dimensional case and to \cite{CaMa} for the case of a thin cylinder. We refer also to \cite{AlFr}, \cite{APR}, \cite{HMC} for the obtention of correctors in other problems related to wave propagation.\par
More general than the wave equation, we consider in this work a Maxwell system. Similarly to the homogenization of (\ref{int1}), one can consider (see e.g. \cite{SaPa} chapter 7 and \cite{Wei})
an open domain $\Om\subset\RR^3$, $\alpha,\beta\in\RR$, $0<|\alpha|+|\beta|$, two sequences  $P_\ep$, $Q_\ep$ of measurable symmetric matrix functions in $\Om$, uniformly elliptic and bounded, two sequences $f_\ep,$ $g_\ep$ in $(0,T)\times \Om$ and two sequences $B_\ep^0$, $D_\ep^0$ in $\Om$ such that there exist $f,\, g,\, B^0$, and $D^0$, satisfying
$$f_\ep\rightharpoonup f\ \hbox{ in }L^2(0,T;L^2(\Om))^3,\quad {\rm div}_xf_\ep\rightarrow {\rm div}_xf\ \hbox{ in }L^2(0,T;H^{-1}(\Om))^3$$
$$g_\ep\rightharpoonup g\ \hbox{ in }L^2(0,T;L^2(\Om))^3,\quad {\rm div}_xg_\ep\rightarrow {\rm div}_xg\ \hbox{ in }L^2(0,T;H^{-1}(\Om))^3$$
$$B^0_\ep\rightharpoonup B^0\ \hbox{ in }L^2(\Om)^3,\quad {\rm div}_xB^0_\ep\rightarrow {\rm div}_xB^0\ \hbox{ in }H^{-1}(\Om)^3$$
$$D^0_\ep\rightharpoonup D^0\ \hbox{ in }L^2(\Om)^3,\quad {\rm div}_xD^0_\ep\rightarrow {\rm div}_xD^0\ \hbox{ in }H^{-1}(\Om)^3.$$
In these conditions, it is proved that the solution $(B_\ep,E_\ep)$ of ($\nu$ denotes the unit normal vector on the boundary, other boundary conditions can be considered)
\beq\label{MSyI}\left\{\ba{l}\dis \partial_t B_\ep+{\rm curl}_x E_\ep=f_\ep\ \hbox{ in }(0,T)\times\Om\\ \ecart\dis
P_\ep\partial_t E_\ep-{\rm curl}_x \big(Q_\ep^{-1}B_\ep\big)=g_\ep\ \hbox{ in }(0,T)\times\Om\\ \ecart\dis
(\alpha B_\ep+\beta E_\ep)\times\nu=0\ \hbox{ on }(0,T)\times \partial\Om\\ \ecart\dis
(B_\ep)_{|t=0}=B^0_\ep\ \  P_\ep E_\ep=D_\ep^0\quad\hbox{ in }\Om,\ea\right.\eeq 
converges weakly in $[L^2(0,T;L^2(\Om))^3]^2$ to the solution $(E,B)$ of the analogous problem by replacing $f_\ep,g_\ep,E_\ep^0$, $D_\ep^0$ by $f,g,E^0$, and $D^0$.
and $P_\ep$, $Q_\ep$ for their $H$-limits $P$, $Q$ (they exist for a subsequence). 
In (\ref{MSyI}), $B_\ep$ and $E_\ep$ represent the magnetic and electric fields respectively. The matrix functions $P_\ep$ and $Q_\ep$ represent the electric permittivity and the magnetic permeability. \par
Despite the previous result, and taking
 $H_\ep^P$ and $H_\ep^Q$ the corrector matrix associated to $P_\ep$ and $Q_\ep$ respectively, it is not true (in general) that 
$$E_\ep-H_\ep^PE\to 0,\ \ B_\ep- Q_\ep H_\ep^QQ_H^{-1}B\to 0\quad\hbox{ in }\, L^1(0,T;L^1(\Om)).$$
It "requires" to assume 
that ${\rm curl}_xE^0_\ep$ and ${\rm curl}_x(Q_\ep^{-1}B_\ep^0)$ converge strongly in $H^{-1}(\Om)^3$ and that the right-hand sides do not oscillate in the time variable.\par 
Related to the results in \cite{BrLe} and \cite{CCMM2} for the wave equation (\ref{int1}), let us assume in this paper that our Maxwell system is posed in the whole of $\RR^3$ and the coefficients $P_\ep$ and $Q_\ep$ have the form $P(x/\ep)$, $Q(x/\ep)$ where $P$ and $Q$ periodic. We notice that we can consider more general coefficients (similar to \cite{CCMM2}) by adding the time dependence but for sake of clarity of the analysis, we choose to consider time-independent coefficients. \par

Assuming that the initial conditions and the right-hand sides of (\ref{MSyI}) have the following form: 
$$B^0\big(x,{x\over\ep}\big),\ D^0\big(x,{x\over\ep}\big),\ f\big(t,x,{t\over\ep},{x\over\ep}\big),\ g\big(t,x,{t\over\ep},{x\over\ep}\big),$$ 
where these functions are  almost periodic in the fast variables $s=t/\ep$ and $y=x/\ep$, we obtain an approximation in the strong topology of $L^2(0,T;L^2(\RR^3))^3$ of $B_\ep$ and $E_\ep$ (corrector result). It consists of adding to the elliptic correctors other oscillating terms which oscillate not only in the space variable but also in the time variable and depend non-locally on the initial data and the second members (see Theorems \ref{thcatwc} and  \ref{Thcorr}). Namely, these terms are written as the sum of infinitely many plane waves in the fast variables $s=t/\ep$ and $y=x/\ep$. Moreover,  the obtained coefficients solve a family of hyperbolic and linear systems depending on the initial conditions and the right-hand sides (see Remark \ref{Rmcae}). \par

Besides of the corrector for the wave system, one of the interests of this work is to show how the Block decomposition method in homogenization can be obtained from the two-scale convergence method for almost periodic functions. For this purpose, we use an extension of the Gelfand transform in $L^2$ to the space of almost periodic functions.
\section{Notation and preliminaries
}\label{SNot}
\begin{itemize}
\item We assume the functions valued in   $\CC$.
\item The ball of center $x\in \RR^N$ and radius $R$ is denoted by $B(x,R)$. 
\item The Lebesgue measure of a measurable set $\om\subset\RR^N$ is denoted by $|\om|$.
\item The real and imaginary part of $z\in \CC$ are respectively denoted by $Re(z)$, $Im(z)$.
The conjugate of $z\in\ZZ$ is denoted by $\overline{z}$.
\item For $\xi,\eta\in \CC^3$, we denote 
$$\xi\cdot\eta:=\sum_{j=1}^3\xi_j\eta_j,$$
where $\xi_j,\eta_j$ denote the components of the vectors $\xi$, $\eta$. The scalar product of $\xi$, $\eta$ is thus given by $\xi\cdot\overline{\eta}$.
\item We denote by $Y$ the unit interval $[-1/2,1/2)$.
\item We denote by $\{e_1,e_2,e_3\}$ the canonical basis in $\RR^3$.
\item We use the index $\sharp$ to mean periodicity with respect to $Y$. For example, $L^2_\sharp(Y^N)$ denotes the spaces of functions in $L^2_{loc}(\RR^N)$ which are periodic, with period $Y^N$.
\item We denote the space of trigonometric functions in $\RR^N$, ${\rm Trig}(\RR^N)$, as the vectorial space spanned by the functions of the form $y\in\RR^N\mapsto e^{i\xi\cdot y}$, $\xi\in \RR^N$. 
\item The closure of ${\rm Trig}(\RR^N)$  in $L^\infty(\RR^N)$ is denoted as $CAP(\RR^N)$. It is the space of almost periodic functions in the Bohr sense. The space of functions in $CAP(\RR^N)$ whose derivatives of order less or equal than $l\in\NN\cup\{+\infty\}$ are also in $CAP(\RR^N)$ is denoted by
$CAP^l(\RR^N)$.
\item The space of almost-periodic functions in the Besicovitch sense, ${\cal B}^2(\RR^N)$ is defined as the space of functions  $u\in L^2_{loc}(\RR^N)$ such that for every $\ep>0$, there exists a trigonometric function $u_\ep\in {\rm Trig}(\RR^N)$ satisfying
$$\limsup_{n\to \infty}{1\over|B(0,n)|}\int_{B(0,n)}|u-u_\ep|^2dy<\ep.$$
\item We say that a function $u\in L^1_{loc}(\RR^N)$ has a mean value $M(u)\in\CC$ if
$$ \lim_{n\to \infty}{1\over |n\om|}\int_{n\om} udy=M(u),\quad\forall\, \om\subset\RR^N, \hbox{ bounded, measurable},\ |\om|>0.$$
\item We recall that for every  $u\in {\cal B}^2(\RR^N)$, there exist $M(u)$ and $M(u^2)$. Identifying  functions $u,v$ in $\RR^N$ such that $M(|u-v|^2)=0$,
the space ${\cal B}^2(\RR^N)$ becomes a Hilbert space with the scalar product
$$(u,v)_{{\cal B}^2(\RR^N)}=M(u\overline v),\quad \forall\, u,v\in {\cal B}^2(\RR^N).$$
Observe that with this identification every function in $L^2(\RR^N)$
is in the class of the null function. Thus,
two functions which are the same as element of ${\cal B}^2(\RR^N)$ can be different in every point. By definition, Trig$(\RR^N)$ and then $CAP^\infty(\RR^N)$ are dense in ${\cal B}^2(\RR^N)$.\par
The  space ${\cal B}^2(\RR^N)$ can be identified with the space of formal series of the form
$$\sum_{\xi\in\RR^N}a_\xi e^{i\xi\cdot y},\quad\hbox{ such that }\ \sum_{\xi\in\RR^N}|a_\xi|^2<\infty,$$
endowed with the scalar product
$$(u,v)=\sum_{\xi\in\RR^N} a_\xi \overline{b_\xi},\quad u=\sum_{\xi\in\RR^N}a_\xi e^{i\xi\cdot y},\ v=\sum_{\xi\in\RR^N}b_\xi e^{i\xi\cdot y}.$$
With this representation, for every $u=\sum_{\xi\in\RR^N}a_\xi e^{i\xi\cdot y}\in{\cal B}^2(\RR^N)$ we have $M(u)=a_0$ and more generally $a_\xi=M(e^{-i\xi\cdot y}u)$, for ever $\xi\in\RR^N$.\par
\item For $m\in\ZZ$, we define ${\cal B}^{m,2}(\RR^N)$ as the space of formal series
$$u=\sum_{\xi\in\RR^N}a_\xi e^{i\xi\cdot y},\ \hbox{ such that }\left\{\ba{ll}\dis  \sum_{\xi\in\RR^N}|a_\xi|^2(1+|\xi|)^{2m}<\infty & \hbox{ if }m\geq 0\\  \ecart\dis\sum_{\xi\in\RR^N}|a_\xi|^2|\xi|^{2m}<\infty& \hbox{ if }m<0.\ea\right.$$
It is a Hilbert space with the scalar product
$$M(u,v)=\left\{\ba{ll}\dis  \sum_{\xi\in\RR^N}a_\xi\overline{b_\xi}(1+|\xi|)^{2m} & \hbox{ if }m\geq 0\\  \ecart\dis\sum_{\xi\in\RR^N}a_\xi\overline{b_\xi}|\xi|^{2m} & \hbox{ if }m<0,\ea\right.\qquad u=\sum_{\xi\in\RR^N}a_\xi e^{i\xi\cdot y},\ v=\sum_{\xi\in\RR^N}b_\xi e^{i\xi\cdot y}.$$
Observe that ${\cal B}^{0,2}(\RR^N)= {\cal B}^{2}(\RR^N)$. For  $u\in{\cal B}^{m,2}(\RR^N)$ we can define the partial derivative $\partial_{y_j} u\in {\cal B}^{m-1,2}(\RR^N)$ by
$$\partial_{y_j}\Big(\sum_{\xi\in\RR^N}a_\xi e^{i\xi\,y}\Big)=i\sum_{\xi\in\RR^N}a_\xi\,\xi_je^{i\xi\cdot y}.$$
\item For $\xi\in Y^N$, we define the operator ${\cal T}_\xi: {\cal B}^2(\RR^N)\to L^2_\sharp(Y^N)$ by
$${\cal T}_\xi u(y)=\sum_{k\in\ZZ^N} M(u(z)e^{-i(\xi+k)\cdot z})e^{ik\cdot y},\quad \forall\, u\in {\cal B}^2(\RR^N).$$
We have
\beq\label{comtp}{\cal T}_\xi(uv)=({\cal T}_\xi u)v,\quad \forall\, u\in {\cal B}^2(\RR^N),\ v\in L^\infty_\sharp(Y^N),\ \xi\in Y^N.\eeq
\beq\label{compde}{\cal T}_\xi (\partial_ju)=i\xi_j{\cal T}_\xi (u)+\partial_j{\cal T}_\xi u,\quad \forall\, u\in {\cal B}^2(\RR^N),\ \xi\in Y^N.\eeq
The operators ${\cal T}_\xi$ allow us to decompose a function in ${\cal B}^2(\RR^N)$ as
\beq\label{desoB} u=\sum_{\xi\in Y^N} e^{i\xi\cdot y}{\cal T}_\xi(u),\eeq
where the functions $\tau_\xi(u)$ are periodic with period $Y^N$.
This is the analogous of the Gelfand transform for functions in $L^2(\RR^N)$ to the case of ${\cal B}^2(\RR^N)$ and agrees with the decomposition of $u$ in Bloch's waves.
\item In the present paper we are mainly interested in almost periodic functions $u=u(s,y)$ in $\RR^4$, with $s$ (the time variable) in $\RR$ and $y$ (the space variable) in $\RR^3$. These functions will be then decomposed as
$$u=\sum_{\lambda\in\RR}\sum_{\xi\in\RR^3}a_{\lambda,\xi}e^{i(\lambda s+\xi\cdot y)}$$
We define the mean values of $u\in{\cal B}^2(\RR^4)$ with respect to $s$ and $y$ by
$$M_s(u)=\sum_{\xi\in\RR}a_{0,\xi}e^{i\xi\cdot y}\in{\cal B}^2(\RR^3),\quad M_y(u)=\sum_{\lambda\in\RR}a_{\lambda,0}e^{i\lambda s}\in{\cal B}^2(\RR).$$
The mean value of $u$ will be denoted by $M_{sy}(u)=a_{0,0}$.\par
We will use the operators ${\cal T}_\xi$ with respect to the space variable. An almost periodic function $u\in{\cal B}^2(\RR^4)$  will be decomposed as
\beq\label{desap} u(s,y)=\sum_{\xi\in Y^3}e^{i \xi\cdot y}{\cal T}_\xi u\eeq
with
$${\cal T}_\xi u:=\sum_{\lambda\in \RR}e^{i\lambda s}{\cal T}_\xi \big(M_s(e^{-i\lambda s}u)\big)\in L^2_\sharp(Y^3; {\cal B}^2(\RR)),\quad
\sum_{\lambda\in\RR\atop\xi\in Y^3}\big\|{\cal T}_\xi\big(M_s(e^{-i\lambda s}u)\big)\big\|_{L^2_\sharp(Y^3)}^2<\infty.$$
\end{itemize}
\section{Position of the problem and main results.}\label{secho}

Our aim in the present section is to get a corrector result for the homogenization of  the Maxwell system
in $\RR^3$ with periodic coefficients. Namely, we are interested in obtaining an approximation in the strong $L^2$-topology of the solutions of 
\beq\label{maxwco}
\left\{\ba{l}\dis \partial_tB_\ep +\curl_x E_\ep=f_\ep\ \hbox{ in }(0,T)\times\RR^3\\ \ecart\dis P\big({x\over\ep}\big)\partial_tE_\ep-\curl_x \big(Q^{-1}\big({x\over \ep}\big)B_\ep\Big)=g_\ep\ \hbox{ in }(0,T)\times\RR^3\\ \ecart\dis (B_\ep)_{|t=0}=B^0_\ep,\
(E_\ep)_{|t=0}=E^0_\ep \ 
\hbox{ in }\RR^3.\ea\right.\eeq
The coefficients $P$ and $Q$ are assumed to satisfy: 
\beq\label{hipcoe}\ba{c}\dis P,\, Q\in L^\infty_\sharp (Y^3)^{3\times 3}\ \hbox{ hermitian}\\ \ecart\dis \exists\,\alpha>0, \ \min\Big\{P\xi\cdot \overline\xi,\ Q\xi\cdot \overline\xi\Big\}\geq \alpha|\xi|^2,\ \forall\, \xi\in\CC^3,\ \hbox{a.e. in }\, Y\ea\eeq
Observe that assuming 
$$f_\ep=0,\quad {\rm div}_xB^0_\ep=0,$$
defining $\rho_\ep$ as the solution of 
$$\left\{\ba{l} \dis \partial_t\rho_\ep={\rm div}_x g_\ep,\ \hbox{ in }(0,T)\times\RR^3\\ \ecart\dis
(\rho_\ep)_{|t=0}={\rm div}_x\Big(P\big({x\over\ep}\big)E^0_\ep\Big), \ \hbox{ in }\RR^3, \ea\right.$$
and applying the divergence (\ref{maxwco}), we get, thanks to initial conditions, 
$${\rm div}_xB_\ep=0,\ \ {\rm div}_x \big(P\big({x\over\ep}\big)E_\ep\big)=\rho_\ep\ \hbox{ in }(0,T)\times \RR^3$$
The magnetic field $B_\ep$ and the electric field $E_\ep$ satisfy the usual Maxwell system
\beq\label{maxwcoU}\left\{\ba{l}\dis {\rm div}_xB_\ep=0\ \hbox{ in }(0,T)\times\RR^3\\ \ecart\dis
 \partial_tB_\ep +\curl_x E_\ep=0\ \hbox{ in }(0,T)\times\RR^3\\ \ecart\dis 
 {\rm div}_x\big(P\big({x\over\ep}\big)E_\ep\big)=\rho_\ep\ \hbox{ in }(0,T)\times \RR^3\\ \ecart\dis
 P\big({x\over\ep}\big)\partial_tE_\ep-\curl_x \Big(Q^{-1}\Big({x\over \ep}\Big)B_\ep\Big)=-j_\ep\ \hbox{ in }(0,T)\times\RR^3\\ \ecart\dis
(E_\ep)_{|t=0}=E^0_\ep,\ \ (B_\ep)_{|t=0}=B^0_\ep\ \hbox{ in }\RR^3,\ea\right.\eeq
where $P(x/\ep)$ is the electric permittivity, $Q(x/\ep)$ the magnetic permeability, 
$\rho_\ep$ and $j_\ep:=-g_\ep$ represent the charge and current densities. \par\medskip
The existence and uniqueness of solutions for (\ref{maxwco}) is based on the following general argument:
\begin{Thm} \label{TexsoM} Assume $P,Q\in L^\infty(\RR^3)^{3\times 3}$ hermitian, satisfying
\beq\label{hipeli}\exists\, \alpha>0,\ \hbox{ such that }\ \min\Big\{P(x)\xi\cdot \overline{\xi},\ Q(x)\xi\cdot \overline{\xi}\Big\}\geq \alpha|\xi|^2,\ \ \forall\, \xi\in\CC^3, \hbox{ a.e. }x\in\RR^3.\eeq
Then, for every $(f,g)\in [L^1(0,T;L^2(\RR^3))^3]^2$ and every $(E^0,B^0)\in [L^2(\RR^3)^3]^2$, there exists a unique solution 
$(E,B)\in [L^\infty(0,T;L^2(\RR^3))^3]^2$ of 
\beq\label{MaxR3}\left\{\ba{l}\dis \partial_tB +\curl_x E=f\ \hbox{ in }(0,T)\times\RR^3\\ \ecart\dis P\partial_tE-\curl_x (Q^{-1}B\Big)=g\ \hbox{ in }(0,T)\times\RR^3\\ \ecart\dis
B_{|t=0}=B^0,\ \ E_{|t=0}=E^0\ \hbox{ in }\RR^3.\ea\right.\eeq
Moreover, 
\beq\label{EnId} {1\over 2}{d\over dt}\int_{\RR^3} \big(Q^{-1}B\cdot \overline{B}+PE\cdot \overline{E})\,dx=\int_{\RR^3}Re(Q^{-1}f\cdot \overline{B}+g\cdot \overline{E}\big)\,dx\ \hbox{ in }(0,T)\eeq
is satisfied 
\end{Thm}\par\medskip\noindent
{\bf Proof.} Using the change of variables $PE=F$, system (\ref{MaxR3}) becomes
\beq\label{MaxR3b}\left\{\ba{l}\dis \partial_tB +\curl_x (P^{-1}F)=f\ \hbox{ in }(0,T)\times\RR^3\\ \ecart\dis \partial_tF-\curl_x (Q^{-1}B\Big)=g\ \hbox{ in }(0,T)\times\RR^3\\ \ecart\dis
B_{|t=0}=B^0,\ \ F_{|t=0}=P^{-1}E^0\ \hbox{ in }\RR^3.\ea\right.\eeq
We define ${\cal A}: D({\cal A})\subset \big[L^2(\RR^3)^3\big]^2\to  \big[L^2(\RR^3)^3\big]^2$,  by
$$D({\cal A})=\big\{(B,F)\in  \big[L^2(\RR^3)^3\big]^2:\ {\curl} (P^{-1}F),{\curl} (Q^{-1}B)\in L^2(\RR^3)^3 \big\}$$
$${\cal A}(B,F)=\big(-{\rm curl}_x (P^{-1}F),{\rm curl}_x (Q^{-1}B)\big),\quad \forall\, (B,F)\in D({\cal A}).$$
We observe that  $D({\cal A})$ is dense in $\big[L^2(\RR^3)^3\big]^2$ and that ${\cal A}$ is closed. Moreover, taking in $\big[L^2(\RR^3)^3\big]^2$ the scalar product
\beq\label{scpr}\big((B_1,F_1),(B_2,F_2)\big)=\int_{\RR^3}\big(Q^{-1}B_1\cdot\overline B_2+P^{-1}F_1\cdot\overline F_2\big)dx,\eeq
we have 
$$\big({\cal A}(B_1,F_1),(B_2,F_2)\big)=\int_{\RR^3}\Big({\rm curl}_x\big(Q^{-1}B_1\big)\cdot\overline{P^{-1}F_2}-
{\rm curl}_x\big(Q^{-1}B_2\big)\cdot\overline{P^{-1}F_1}\Big)dx$$
 for every $(B_1,F_1) and (B_2,F_2)\in D({\cal A})$.
This proves that ${\cal A}$ is anti-hermitic, which combined with the density of $D({\cal A})$ and ${\cal A}$ closed allows us to apply the Lax-Milgram theorem to deduce that   ${\cal A}$ is maximal monotone. Therefore the Hille-Yoshida theorem provides the existence of solution for  (\ref{MaxR3b}) and proves (\ref{EnId}) in the new variables $B$, $F$. Returning to the variables $B$, $F$ we have the thesis of Theorem \ref{TexsoM}. \cqfd
\begin{Rem}\label{EpMa} Let 
$$H(t)=\int_{\RR^3} \big(Q^{-1}B\cdot \overline{B}+PE\cdot \overline{E})\,dx.$$
Equation (\ref{EnId}) implies
$${1\over 2}{dH\over dt}\leq \left(\int_{\RR^3}\big(Q^{-1}f\cdot \overline{f}+P^{-1}g\cdot \overline{g}\big)dx\right)^{1\over 2}H^{1\over 2},$$
and then
$${dH^{1\over 2}\over dt}\leq \left(\int_{\RR^3}\big(Q^{-1}f\cdot \overline{f}+P^{-1}g\cdot \overline{g}\big)dx\right)^{1\over 2}.$$
Integrating this inequality we deduce the a priori estimate
\beq\label{apre} \ba{l}\dis\left(\int_{\RR^3} \big(Q^{-1}B\cdot \overline{B}+PE\cdot \overline{E})\,dx\right)^{1\over 2}\leq 
 \left(\int_{\RR^3} \big(Q^{-1}B^0\cdot B^0+PE^0\cdot E^0)\,dx\right)^{1\over 2}\\ \ecart \qquad\qquad\qquad\qquad\qquad\qquad\dis+\int_0^t\left(\int_{\RR^3}\big(Q^{-1}f\cdot \overline{f}+P^{-1}g\cdot \overline{g}\big)dx\right)^{1\over 2}ds\ \hbox{ in }(0,T).\ea
\eeq
\end{Rem}\par\medskip
 Taking into account Theorem \ref{TexsoM}, let us  assume in  (\ref{maxwco}) 
\beq\label{hipco} f_\ep,\,g_\ep\ \hbox{ bounded in } L^2(0,T;L^2(\RR^3))^3,\quad E^0_\ep,\,B^0_\ep\ \hbox{ bounded in }L^2(\RR^3)^3.\eeq
Moreover, for the Maxwell system, which is the case we are mainly interested in, we have that ${\rm div}_x B_\ep=0$ in $(0,T)\times \RR^3$ and that ${\rm div}_x(P_\ep E_\ep)$ is the charge which must be finite.  In particular, the divergence of $B_\ep$ and $P_\ep E_	\ep$ must be  bounded. To this aim, let us also assume
\beq\label{consd} \left\{\ba{l}\dis {\rm div}_xf_\ep,\ {\rm div}_xg_\ep\ \hbox{ bounded in }L^2(0,T;L^2(\RR^3))\\ \ecart\dis {\rm div}_x\big(P\big({x\over\ep}\big)E^0_\ep\big),\ {\rm div}_xB^0_\ep\ \hbox{ bounded in }L^2(\RR^3).\ea\right.\eeq\par
From Theorem \ref{TexsoM} and taking the divergence in the two differential equations in (\ref{maxwco}), 
we conclude that (\ref{maxwco}) has a unique solution $(B_\ep,E_\ep)$ satisfying 
\beq\label{estpri} \ba{l}\dis\|B_\ep\|_{L^\infty(0,T;L^2(\RR^3))^3}+\|E_\ep\|_{L^\infty(0,T;L^2(\RR^3))^3}\\ \ecart\dis+
\|{\rm div}_xB_\ep\|_{H^1(0,T;L^2(\RR^3))}+
\Big\|{\rm div}_x\big(P\big({x\over\ep}\big)E_\ep\big)\Big\|_{H^1(0,T;L^2(\RR^3))}\leq C,\ea\eeq
with  $C>0$ independent of $\ep$. This gives the existence of subsequences  converging weakly-$\ast$ in $L^\infty(0,T;L^2(\RR^3))^3$. For our purpose, we are more interested in the existence of two-scale limits. We recall the corresponding definition.
\begin{Def} Let $\Om$ be an open set in $\RR^N$. We say that a sequence $u_\ep$ in $L^2(\Om)$ two-scale converges to a function $ \hat u\in L^2(\Om;{\cal B}^2(\RR^N))$ and we write
$$u_\ep\stackrel{2e}\rightharpoonup \hat u\ \hbox{ in }\Om,$$
if for every $\psi\in C^0_c(\Om;CAP(\RR^N))$, we have
\beq\label{defco2e} \lim_{\ep\to 0}\into u_\ep \psi\big(x,{x\over \ep}\big)dx=\into M(\hat u\psi)dx.\eeq
\end{Def}
The following theorem has been proved in \cite{CaGa}, \cite{CaGa2} (see also \cite{All}, \cite{Ngu} for the periodic case).
\begin{Thm} \label{Thco2e} For every bounded sequence in $L^2(\Om)$, there exists a subsequence of $(u_\ep)_\ep$, still denoted by $u_\ep$, and $\hat u\in L^2(\Om;{\cal B}^2(\RR^N))$ such that $u_\ep$ two-scale converges to $\hat u$. Moreover $u_\ep$ converges weakly in $L^2(\Om)$ to $u:=M(\hat u)$.
\end{Thm}
We notice that the idea behind the definition of  the two-scale convergence is to obtain an approximation of $u_\ep$ of the type
\beq\label{decor} u_\ep(x)\sim u\Big(x,{x\over\ep}\Big).\eeq
\begin{Def} \label{st2s} We say that the sequence $u_\ep\in L^2(\Om)$  two-scale strongly converges to $\hat u\in L^2(\Om;{\cal B}^2(\RR^N))$ and we write
$$u_\ep\stackrel{2e}\rightarrow \hat u\ \hbox{ in }\Om,$$
if for every $\hat v \in L^2(\Om;{\cal B}^2(\RR^N))$ and every 
$v_\ep\in L^2(\Om)$ which two-scale converges to $\hat v$, we have
\beq\label{limddtw}\lim_{\ep\to 0}\into u_\ep v_\ep\,dx=\into M(\hat u\hat v)\,dx.\eeq
\end{Def}
\begin{Rem} \label{stcyc} One can check in analogy with the periodic case (see \cite{All}) that $u_\ep$ two-scale strongly converges to $\hat u$ if and only if
$$\lim_{\ep\to 0} \into |u_\ep|^2dx=\into M(|\hat u|^2)dx.$$
\end{Rem}
Related to (\ref{decor}), another characterization of the strong two-scale convergence which will provide us with the corrector result for  (\ref{maxwco}) states as follows
\begin{Pro} \label{Procf}The sequence $u_\ep$ two-scale strongly converges to $\hat u$ if and only if for every sequence $\hat u^n\in L^2(\Om; CAP(\RR^N))$ which converges to $\hat u\in L^2(\Om;{\cal B}^2(\RR^N))$, we have
\beq\label{confucg} \lim_{n\to\infty}\limsup_{\ep\to 0}\into\Big|u_\ep-\hat u^n\big(x,{x\over\ep}\big)\Big|^2dx=0.\eeq
\end{Pro}
\begin{Rem}\label{re2ef} If in Proposition \ref{Procf} the function $\hat u$ belongs to $L^2(\Om; CAP(\RR^N))$, then we can take
$\hat u^n=\hat u$ for every $n$, to get that
$$u_\ep-\hat u\Big(x,{x\over\ep}\Big)\to 0\ \hbox{ in }L^2(\Om).$$
However, for an arbitrary function $\hat u\in L^2(\Om;{\cal B}^2(\RR^N))$, the function $x\to \hat u(x,x/\ep)$ is not even well defined. Recall that the elements of ${\cal B}^2(\RR^N)$ are not true functions but just class of functions where two functions in the same class can be different in every point. Due to this difficulty, we shall use a double approximation in (\ref{confucg}).
\end{Rem}
Taking into account (\ref{hipco}) and Theorem \ref{Thco2e}  we deduce that, up to extraction of subsequence of $\ep$, still denoted by $\ep$, there exist $\hat f,\,\hat g\in L^2(0,T;L^2(\RR^3;{\cal B}^2(\RR^4))^3$, $\hat B^0,	\hat E^0\in L^2(\RR^3;{\cal B}^2(\RR^3))^3$ such that
\beq\label{conv2eda} \left\{\ba{l}\dis f_\ep\stackrel{2e}\rightharpoonup \hat f,\ \ g_\ep\stackrel{2e}\rightharpoonup \hat g\ \ \hbox{ in }(0,T)\times \RR^3,\\ \ecart\dis
B^0_\ep\stackrel{2e}\rightharpoonup \hat B^0,\ \ E_\ep\stackrel{2e}\rightharpoonup \hat E^0\ \ \hbox{ in }\RR^3.\ea\right.\eeq
Taking such subsequence, it is not restrictive in the following to assume that (\ref{conv2eda}) holds. We recall that
taking
\beq\label{limded}  f=M_{sy}(\hat f),\quad g=M_{sy}(\hat g),\quad B^0=M_y(\hat B^0),\quad (PE)^0=M_y(P\hat E^0),\eeq
we also have
\beq\label{convdd}\begin{array}{c}
f_\ep\rightharpoonup f,\ \ g_\ep\rightharpoonup g\quad\hbox{ in }L^2(0,T;L^2(\RR^3))^3\\ \ecart\dis
P\big({x\over \ep}\big)E^0_\ep\rightharpoonup (PE)^0,\ \ B^0_\ep\rightharpoonup B^0\quad\hbox{ in }L^2(\RR^3)^3.\ea\eeq
\par
Moreover, we also prove, using (\ref{consd}), that  
\begin{Pro}\label{divnli} The functions $\hat f,\hat g, \hat B^0$, and $\hat E^0$ in (\ref{conv2eda}) satisfy
\beq\label{condivd} {\rm div}_y\hat f={\rm div}_y\hat g=0\ \hbox{ in }(0,T)\times\RR^3\times\RR^4,\quad
{\rm div}_y\hat B^0={\rm div}_y\big(P\hat E^0\big)=0\ \hbox{ in }\RR^3\times\RR^3.\eeq
\end{Pro}\par\noindent
{\bf Proof.} By (\ref{consd}), for $\varphi\in C^\infty_c(\RR^3)$, $\lambda\in\RR$, and $\xi\in\RR^3$ we have
$$\ba{l}\dis 0=\lim_{\ep\to 0}\ep \int_0^T\int_{\RR^3} {\rm div}_x f_\ep\,\varphi(x)e^{-i(\lambda {t\over\ep}+\xi\cdot {x\over\ep})}dxdt\\ \ecart\dis
=-\lim_{\ep\to 0}\ep\int_0^T\int_{\RR^3}f_\ep\cdot \nabla\varphi \, e^{-i(\lambda {t\over\ep}+\xi\cdot {x\over\ep})}\,dxdt+\lim_{\ep\to 0}\int_0^T\int_{\RR^3}i\xi\cdot f_\ep\, \varphi e^{-i(\lambda {t\over\ep}+\xi\cdot {x\over\ep})}\,dxdt\\ \ecart\dis
=i\int_0^T\int_{\RR^3}M_{sy}\big(\hat f\cdot\xi\,e^{-i(\lambda s+\xi\cdot y)}\big) \varphi \,dxdt.
\ea$$
By the definition of derivative in the space of almost periodic functions given above, this is equivalent to
${\rm div}_y \hat f=0$. Using the same idea with $g_\ep$, $B^0_\ep$ and $P(x/\ep)E^0_\ep$, we deduce
(\ref{condivd}).\hfill \cqfd\par\
\par
In order to characterize the two-scale limits of the solution $(B_\ep,E_\ep)$ of (\ref{maxwcoU}) let us introduce  the following space 
\begin{Def}\label{decaM} We define ${\cal M}$ as the space of functions $(U,V)\in[{\cal B}^2(\RR^3)^3]^2$ which satisfy the Maxwell system
\beq\label{MSmv}\left\{\ba{l} {\rm div}_yU=0\ \hbox{ in }\RR^4\\ \ecart\dis
\partial_sU+{\rm curl}_yV=0\ \hbox{ in }\RR^4\\ \ecart\dis {\rm div}_y\big(PV)=0\ \hbox{ in }\RR^4\\ \ecart\dis
P\partial_sV-{\rm curl}_y\big(Q^{-1}U\big)=0\ \hbox{ in }\RR^4.\ea\right.
\eeq
\end{Def}
\begin{Rem} The second and fourth equations in (\ref{MSmv}) imply that the first and third equations are equivalent to
\beq\label{sideM}{\rm div}_yM_s(U)={\rm div}_y(PM_s(V))=0\ \hbox{ in }\RR^3.\eeq
\end{Rem}
\par\medskip
Using the decomposition (\ref{desap}) of an element of ${\cal B}^2(\RR^4)$, and taking into account
(\ref{comtp}) and (\ref{compde}), we have that a pair $(U,V)$ belongs to ${\cal M}$ if and only if for every $\xi\in Y^3$, the pair $({\cal T}_\xi U,{\cal T}_\xi V)\in \big[{\cal B}^2(\RR;L^2_\sharp(Y^3))^3\big]^2$ satisfies
\beq\label{MSmvB}\left\{\ba{l} i\xi\cdot{\cal T}_\xi U+{\rm div}_y{\cal T}_\xi U=0\ \hbox{ in }\RR^4\\ \ecart\dis
\partial_s{\cal T}_\xi U+i\xi\times {\cal T}_\xi V+{\rm curl}_y{\cal T}_\xi V=0\ \hbox{ in }\RR^4
\\ \ecart\dis i\xi\cdot P{\cal T}_\xi V+{\rm div}_y\big(P{\cal T}_\xi V\big)=0\ \hbox{ in }\RR^4\\ \ecart\dis
P\partial_s{\cal T}_\xi V-i\xi\times \big(Q^{-1}{\cal T}_\xi U\big)-
{\rm curl}_y\big(Q^{-1}{\cal T}_\xi U\big)=0\ \hbox{ in }\RR^4.\ea\right.
\eeq
This result allows us to get a spectral decomposition of the elements of ${\cal M}$, which is given by 
\begin{Pro}\label{ThcAu}  Assume  $P$ and $Q$ satisfy (\ref{hipcoe}) and consider $\xi\in Y^3$,
then we have
\begin{enumerate}
\item The set of eigenvalues $\lambda\in\RR$ such that the system
\beq\label{pbaut}\left\{\ba{l} U,V\in L^2_\sharp(Y^3)^3\\ \ecart\dis
i\xi\cdot U+{\rm div}_yU=0\ \hbox{ in }\RR^3\\ \ecart\dis
i\lambda U+i\xi\times V+{\rm curl}_y V=0\ \hbox{ in }\RR^3
\\ \ecart\dis i\xi\cdot PV+{\rm div}_y\big(PV\big)=0\ \hbox{ in }\RR^3\\ \ecart\dis
i\lambda PV-i\xi\times Q^{-1}U-
{\rm curl}_y\big(Q^{-1}U\big)=0\ \hbox{ in }\RR^3,\ea\right.
\eeq
has a non-null solution, is given by a sequence $\lambda^\xi_j\in\RR$, with $|\lambda^\xi_j|$ converging to infinity. Moreover, for every $\lambda=\lambda^\xi_j$ the  space ${\cal V}^\xi_{j}$ of solutions of (\ref{pbaut}) has  finite-dimension $m^\xi_j\geq 1$.
\item Taking a basis $(\Phi^\xi_{jl},\Psi^\xi_{jl})$ of ${\cal V}^\xi_{j}$, orthonormal with respect to the scalar product
\beq\label{proeL2} \big((U_1,V_1),(U_2,V_2)\big)_{[L^2_\sharp(Y^3)^3]^2 }=\int_{Y^3} \big(Q^{-1}U_1\cdot \overline{U_2}+
PV_1\cdot \overline{V_2}\big)dy,\eeq
we have that the family 
$$\big\{(\Phi^\xi_{jl},\Psi^\xi_{jl})\big\},\quad j\in\NN,\ 1\leq l\leq m^\xi_j,$$
is an orthonormal basis of
\beq\label{defeWx}W_\xi:=\Big\{(U,V)\in \big[L^2_\sharp(Y^3)^3\big]^2:\ 
i\xi\cdot U+{\rm div}_yU=i\xi\cdot PV+{\rm div}_y(PV)=0\Big\}.\eeq
\item The family 
$$\big\{e^{i\xi\cdot y}(\Phi^\xi_{jl},\Psi^\xi_{jl})\big\},\quad \xi\in Y^3,\ j\in\NN,\ 1\leq l\leq m^\xi_j,$$
is an orthonormal basis of
\beq\label{defeW} \Big\{(U,V)\in \big[{\cal B}^2(\RR^3)^3\big]^2:\ 
{\rm div}_yU={\rm div}_y(PV)=0\ \hbox{ in }\RR^3\Big\}.\eeq\par
\item A pair $(\hat U,\hat V)$ belongs to ${\cal M}$ if and only if it can be decomposed as
\beq\label{descM} (\hat U,\hat V)=\sum_{\xi\in Y^3}\sum_{j=1}^\infty e^{i(\lambda^\xi_j s+\xi\cdot y)} \sum_{l=1}^{m^\xi_j}a^\xi_{jl}\big(\Phi^\xi_{jl}(y),\Psi^\xi_{jl}(y)\big),\qquad
\sum_{\xi\in Y^3}\sum_{j=1}^\infty\sum_{l=1}^{m^\xi_j}|a^\xi_{jl}|^2<\infty.\eeq\par
\item The value $\lambda=0$  belongs to the set of eigenvalues $\lambda_j^\xi$ associated to $\xi$ if and only if $\xi$ is the null vector. In this case the space of solutions of  (\ref{pbaut}) is characterized as the set of functions
\beq\label{solaa0} U=Q(\eta+\nabla_yW_\eta),\quad V=\zeta+\nabla_yZ_\zeta,\eeq
with
\beq\label{MUMV1} \left\{\ba{l}\dis \eta,\zeta\in\RR^3,\ W_\eta,Z_\eta\in H^1_\sharp(Y^3)\\ \ecart\dis -{\rm div}_y\big(Q(\eta+\nabla_yW_\eta)\big)=-{\rm div}_y\big(P(\zeta+\nabla_yZ_\zeta)\big)=0\hbox{ in }\RR^3.\ea\right.\eeq
\end{enumerate}
\end{Pro}\par\medskip\noindent
{\bf Proof.} In order to prove the first and second parts in Proposition \ref{ThcAu} we need to get a spectral decomposition of the operator ${\cal R}:D({\cal R})\subset W_\xi\to W_\xi$ defined by
\beq\label{decR}D({\cal R})=\big\{(U,V)\in W_\xi:\ {\rm curl}_y(Q^{-1}U),{\rm curl}_yV\in L^2_\sharp(Y^3)^3\big\}.\eeq
\beq\label{defOR} {\cal R}(U,V)=\big(i\xi\times V+{\rm curl}_yV, -i\xi\times Q^{-1}U-{\rm curl}_y(Q^{-1}U)\big),\quad \forall\, (U,V)\in D({\cal R}).\eeq
We observe that $D(R)$ is dense in $W_\xi$ and  ${\cal R}$ is anti-hermitic taking $W_\xi$ endowed with the scalar product defined by (\ref{proeL2}). Thus, such spectral decomposition exists but it could be given by an integral decomposition. Using the Fredholm theory, let us prove that there is only a sequence $\lambda_j^\xi$ of singular values of ${\cal B}$ which tends to infinity, and that every of them is an eigenvalue with finite dimension.\par\medskip
For $\xi\in Y^3$, we define ${\cal A}: W_\xi\to W_\xi$ by ${\cal A}(U,V)=(\hat U,\hat V),$
with $(\hat U,\hat V)$ the unique solution of 
$$\left\{\ba{l}\dis (\hat U,\hat V)\in  W_\xi\\ \ecart\dis
 \hat U+i\xi\times \hat V+{\rm curl}_y\hat V=U\ \hbox{ in }\RR^3\\ \ecart\dis
P\hat V-i\xi\times Q^{-1}\hat U-{\rm curl}_y(Q^{-1}\hat U)=PV\ \hbox{ in }\RR^3.\ea\right.$$
The existence and uniqueness of solution for this problem follows writing the equations as
$$\left\{\ba{l} \hat U+i\xi\times P^{-1}\hat F+{\rm curl}_y(P^{-1}\hat F)=U\ \hbox{ in }\RR^3\\ \ecart\dis
\hat F-i\xi\times Q^{-1}\hat U-{\rm curl}_y(Q^{-1}\hat U)=F\ \hbox{ in }\RR^3\ea\right.$$
with $\hat F=P\hat V$, $F=PV$ and observing that the operator
${\cal B}:D(B)\subset \big[L^2(Y^3)^3\big]^2\to \big[L^2(Y^3)^3\big]^2$ defined by
$$D({\cal B})=\Big\{(U,F)\in \big[L^2(Y^3)^3\big]^2:\ {\rm curl}_y(Q^{-1}U),{\rm curl}_y(P^{-1}F)\in L^2(Y^3)^3\Big\},$$
$${\cal B}(F,U)=\Big(i\xi\times P^{-1}F+{\rm curl}_y(P^{-1}F),-i\xi\times Q^{-1} U-{\rm curl}_y(Q^{-1}U)\Big),\quad \forall\,(F,U)\in D({\cal B}),$$
has dense domain and is closed and anti-hermitic with respect to the scalar product in $\big[L^2(Y^3)^3\big]^2$ defined by
$$\big((U_1,F_1),(U_2,F_2)\big)=\int_{Y^3} \big(Q^{-1}U_1\cdot \overline U_2+P^{-1}F_1\cdot F_2\big).$$
This implies that ${\cal B}$ is maximal monotone and then the existence of $(\hat U,\hat V)$.\par
Let us prove that ${\cal A}$ is compact. We consider a bounded sequence $(U_n,V_n)$ in $W_\xi$. Using that $(\hat U_n,\hat V_n):={\cal A}(U_n,V_n)$ satisfies
$$\int_{Y^3}\big(Q^{-1}\hat U_n\cdot\overline{\hat U_n}+P\hat V_n\cdot\overline{\hat V_n}\big)dy=\int_{Y^3}\big(Q^{-1}\hat U\cdot\overline{\hat U}+P\hat V\cdot\overline{\hat V}\big)dy$$
we get that $(\hat U_n,\hat V_n)$ is also bounded in $W_\xi$. From definition of ${\cal A}$, we also have that
${\rm curl}_y\hat V_n$, ${\rm curl}_y(Q^{-1}\hat U_n)$ are bounded in $L^2(Y^3)^3$. On the other hand, using a Hodge decomposition for $Q^{-1}\hat U_n$ and $\hat V_n$, we have the existence of  $C>0$,
$\hat W_n,\hat Z_n\in H^1_\sharp(Y^3)^3$ and $\hat w_n,\hat z_n\in H^1_\sharp(Y^3)$ such that
\beq\label{dthco1}Q^{-1}\hat U_n=\hat W_n+\nabla_y \hat w_n,\quad \hat V_n=\hat Z_n+\nabla_y \hat z_n,\quad {\rm div}_y\hat W_n={\rm div}_y(Q^{-1}\hat Z_n)=0,\eeq
with
$$\|\hat W_n\|_{H^1(Y^3)^3}+\|\hat w_n\|_{H^1(Y^3)}\leq C\Big(\|Q^{-1}\hat U_n\|_{L^2(Y^3)^3}+\|{\rm curl}_y\big(Q^{-1}\hat U_n\|_{L^2(Y^3)^3}\Big),$$
$$\|\hat Z_n\|_{H^1(Y^3)^3}+\|\hat z_n\|_{H^1(Y^3))}\leq C\Big(\|\hat V_n\|_{L^2(Y^3)^3}+\|{\rm curl}_y\hat V_n\|_{L^2(Y^3)^3}\Big).$$\par
Using (\ref{dthco1}) and $(\hat U_n,\hat V_n)\in W_\xi$ we get
$${\rm div}_y(Q\nabla_y \hat w_n)=-i\xi\cdot \hat U_n-{\rm div}_y(Q\hat W_n),\quad {\rm div}_y(P\nabla_y \hat z_n)=-i\xi\cdot P\hat V_n-{\rm div}_y(P\hat Z_n),$$
By Rellich-Kondrachov's compactness theorem $\hat W_n$, $\hat Z_n$ are compact in $L^2_\sharp(Y^3)$ and then
the right-hand sides of these equalities are compact in $H^{-1}_\sharp(Y^3)$. Thus,  by Lax-Milgram's theorem  $\hat w_n$ and $\hat z_n$  are compact in $H^1_\sharp(Y^3)$. Returning to (\ref{dthco1}) we conclude that $\hat U_n$, $\hat V_n$ are compact in $L^2(Y^3)^3$.  This proves the compactness of ${\cal A}$.\par\medskip
In order to prove that the spectral values of ${\cal R}$ are a sequence of eigenvalues with finite dimension which tends to infinity it is enough to observe that
$${\cal R}={\cal P}(A^{-1}-I),$$
with ${\cal P}:\big[L^2(\RR^3)^3\big]^2\to \big[L^2(\RR^3)^3\big]^2$ defined by
$${\cal P}(U,V)=(U,PV),\quad \forall\, (U,V)\in \big[L^2(\RR^3)^3\big]^2.$$
Thus, the result follows from the Fredholm theory.\par\medskip
The third part of Proposition \ref{ThcAu} is a simple consequence of the second one and decomposition 
(\ref{desoB}) of a function in ${\cal B}^2(\RR^3)$. Analogously, the fourth part follows from decomposition (\ref{desap}) for functions in ${\cal B}^2(\RR^4)$.\par\medskip
Finally, let us now prove the fifth part of the proposition. We have that $\lambda=0$ is an eigenvalue associated to $\xi$ if and only if there exists a non-null vector $(U,V)\in L^2_\sharp(Y^3)^3$ solution of 
$$\left\{\ba{l} 
i\xi\cdot U+{\rm div}_yU=i\xi\cdot PV+{\rm div}_y\big(PV\big)=0\ \hbox{ in }\RR^3\\ \ecart\dis
i\xi\times V+{\rm curl}_y V=i\xi\times Q^{-1}U+
{\rm curl}_y\big(Q^{-1}U\big)=0\ \hbox{ in }\RR^3,\ea\right.
$$
or equivalently
\beq\label{conxi0}\left\{\ba{l}\dis {\rm curl}_y\big(e^{i\xi\cdot y}Q^{-1}U\big)={\rm curl}_y\big(e^{i\xi\cdot y}V\big)=0\ \hbox{ in }\RR^3,\\ \ecart\dis {\rm div}_y\big(e^{i\xi\cdot y}U\big)={\rm div}_y\big(e^{i\xi\cdot y}PV\big)=0\ \hbox{ in }\RR^3.\ea\right.\eeq
If $\xi$ is not the null vector, this implies the existence of $W,Z\in H^1_\sharp(Y^3)$ such that
$$e^{i\xi\cdot y}U=Q\nabla_y\big(e^{i\xi\cdot y} W\big),\quad {\rm div}_y\Big(Q\nabla_y\big(e^{i\xi\cdot y} W\big)\big)=0\ \hbox{ in }\RR^3$$ 
$$e^{i\xi\cdot y} V=\nabla_y\big(e^{i\xi\cdot y} Z\big),\quad {\rm div}_y\Big(P\nabla_y\big(e^{i\xi\cdot y}Z\big)\big)=0\ \hbox{ in }\RR^3,$$ 
but this is only possible if $W$ and $Z$ are the null function what leads to the contradiction $U=V=0$.\par
If $\xi=0$ then (\ref{conxi0}) proves  the existence of $\eta,\zeta\in\RR^3$ such that (\ref{solaa0}) and (\ref{MUMV1}) hold. \cqfd
\begin{Rem} \label{caRe} If $P$ and $Q$ are real functions, which is the case for Maxwell system, then
the set $\{\lambda^\xi_j\}$ is stable for the change of sign and the space of solutions of (\ref{pbaut}) for $-\lambda^\xi_j$ agrees with the conjugated of the space corresponding to $\lambda^\xi_j$.
\end{Rem}
\begin{Rem} The functions $W_\eta$ and $Z_\zeta$ in (\ref{MUMV1}) are the usual corrector terms associated to the functions $x\mapsto \eta\cdot x$ and $ x\mapsto \zeta\cdot x$ relative to the homogenization of  the elliptic operators 
$-{\rm div }\big(Q\big({x\over\ep}\big)\nabla u_\ep\big)$
and 
$-{\rm div }\big(P\big({x\over\ep}\big)\nabla u_\ep\big)=G$
respectively (see e.g. \cite{All}, \cite{BLP},  \cite{Ngu}, \cite{Tar}). \end{Rem}
We are now in position of giving the following theorem which provides the two-scale limits of the solutions $(B_\ep,E_\ep)$ of (\ref{maxwcoU}).
\begin{Thm} \label{thcatwc} Assume $P$, $Q$ satisfy (\ref{hipcoe}). Then, for every $f_\ep,g_\ep\in L^2(0,T;L^2(\RR^3))^3$, $B^\ep_0,E^\ep_0\in L^2(\RR^3)^3$ which satisfy (\ref{consd}) and are such that there exist $\hat f,\hat g\in L^2(0,T;L^2(\RR^3;{\cal B}^2(\RR^4))^3$, $\hat B^0,	\hat E^0\in L^2(\RR^3;{\cal B}^2(\RR^3))^3$ satisfying (\ref{conv2eda}), we have that the unique solution $(B_\ep,E_\ep)$ of (\ref{maxwco}) satisfies
\beq\label{conv2eso} B_\ep\stackrel{2e} \rightharpoonup \hat B,\ \
E_\ep\stackrel{2e} \rightharpoonup \hat E\quad \hbox{in }(0,T)\times\RR^3,\eeq
where, for ${\cal M}$ given by Definition \ref{decaM}, the functions
$\hat B,\hat E\in L^\infty(0,T; L^2(\RR^3;{\cal B}^2(\RR^4)))^3$ are characterized by
\beq\label{pbli}\left\{\ba{l}\dis (\hat B,\hat E)\in L^\infty(0,T;L^2(\RR^3;{\cal M}))\\ \ecart\dis
\partial_tM_{sy}\big(Q^{-1}\hat B\cdot \overline U+P\hat E\cdot\overline V\big)+{\rm div}_xM_{sy}\big(\hat E\times \overline {Q^{-1}U}-Q^{-1}\hat B\times \overline V\big)\\ \ecart\dis \qquad=M_{sy}\big(Q^{-1}\hat f\cdot \overline U+\hat g\cdot \overline V\big)\ \hbox{ in }(0,T)\times \RR^3\\ \ecart\dis
M_{sy}\big(Q^{-1}\hat B\cdot\overline U+P\hat E\cdot\overline V\big)_{|t=0}=M_{y}\big(Q^{-1}\hat B^0\cdot\overline U_{|s=0}+P\hat E^0\cdot\overline V_{|s=0}\big),\ \forall\, (U,V)\in {\cal M}.
\ea\right.\eeq
\end{Thm}\par\noindent
{\bf Proof.} By (\ref{estpri}), there exist $\hat B,\hat E\in L^2(0,T; L^2(\RR^3;{\cal B}^2(\RR^4)))^3$ such that for a subsequence of $\ep$, still denoted by $\ep$, (\ref{conv2eso}) holds. Taking into account that $B_\ep$, $E_\ep$ are bounded in $L^\infty(0,T;L^2(\RR^3))^3$, it is simple to check that  $\hat B,\hat E$ belong to $L^\infty(0,T; L^2(\RR^3;{\cal B}^2(\RR^4)))^3$. \par
For $U\in C^\infty_c(0,T; C^\infty_c(\RR^3;CAP^\infty(\RR^4)))^3$,  we take 
$\ep U(t,x,t/\ep,x/\ep)$
as test function in the first equation in (\ref{maxwco}). This provides
$$\ba{l}\dis 
-\int_0^T\hskip-3pt\int_{\RR^3} B_\ep\cdot\Big(\ep\,\partial_tU\big(t,x,{t\over\ep},{x\over\ep}\big)+\partial_sU\big(t,x,{t\over\ep},{x\over\ep}\big)\Big)dxdt\\ \ecart\dis
+\int_0^T\hskip-3pt\int_{\RR^3} E_\ep\cdot\Big(\ep\,{\rm curl}_xU\big(t,x,{t\over\ep},{x\over\ep}\big)+{\rm curl}_yU\big(t,x,{t\over\ep},{x\over\ep}\big)\Big)dxdt \\ \ecart\dis =\ep\int_0^T\hskip-3pt\int_{\RR^3} f_\ep\cdot U\big(t,x,{t\over\ep},{x\over\ep}\big)\Big)dxdt.
\ea$$
Passing to the limit in this equality thanks to $U_\ep, E_\ep$ bounded in $L^\infty(0,T;L^2(\RR^3))^3$ and (\ref{conv2eso}) we get
$$
\int_0^T\hskip-3pt\int_{\RR^3}M_{sy}\big(- \hat B\cdot\partial_sU+\hat E\cdot {\rm curl}_yU\big)dxdt=0,\quad \forall\,U\in C^\infty_c(0,T; C^\infty_c(\RR^3;CAP^\infty(\RR^4)))^3.$$
This proves that $(\hat B,\hat E)$ satisfies
$$\partial_s\hat B+{\rm curl}_y\hat E=0.$$
A similar reasoning, using the second equation in (\ref{maxwco}) and  ${\rm div}_xB_\ep$, ${\rm div}_y(P(x/\ep)E_\ep)$ bounded in $L^\infty(0,T;L^2(\RR^3))$ also implies
$$P\partial_s\hat E-{\rm curl}_y(Q^{-1}\hat B)=0,\quad {\rm div}_y\hat B=0,\quad {\rm div}_y(P\hat E)=0,$$
Thus,
\beq\label{limpeM} \hat B,\hat E\in {\cal M}\ \hbox{ a.e in }(0,T)\times \RR^3,\eeq
with ${\cal M}$ given by Definition \ref{decaM}.\par
Now, for $\xi\in Y^N$, $j\in\NN$ and $(U^\xi_j,V^\xi_j)\in C^\infty_c([0,T);C^\infty_c(\RR^3;{\cal V}^\xi_j))$ (see the statement of  Proposition \ref{ThcAu} for the definition of $\lambda^\xi_j$, ${\cal V}^\xi_j$) we define
$$U(t,x,s,y)=e^{i(\lambda^\xi_j s+\xi\cdot y)}U^\xi_j(t,x,y),\quad
V(t,x,s,y)=e^{i(\lambda^\xi_js+\xi\cdot y)}V^\xi_j(t,x,y).$$
Taking
$$\overline {Q^{-1}}\big({x\over\ep}\big)\overline{U}\big(t,x,{t\over \ep},{x\over\ep}\big),$$
as test function in the first equation in (\ref{maxwco}), we get
\beq\label{e1dtp}\ba{l}\dis -\int_{\RR^3}Q^{-1}\big({x\over\ep}\big)B_\ep^0\cdot \overline U\big(0,x,0,{x\over\ep}\big)dx\\ \ecart\dis-\int_0^T\hskip-3pt\int_{\RR^3}Q^{-1}\big({x\over\ep}\big)B_\ep\cdot \Big( \partial_t\overline U\big(t,x,{t\over \ep},{x\over\ep}\big)+{1\over\ep}\partial_s\overline U\big(t,x,{t\over \ep},{x\over\ep}\big)\Big)dxdt\\ \ecart\dis+\int_0^T\hskip-3pt\int_{\RR^3}E_\ep\cdot \Big( {\rm curl}_x\Big(\overline{Q^{-1}}\big({x\over\ep}\big)\overline U\big(t,x,{t\over \ep},{x\over\ep}\big)\Big)+{1\over\ep}{\rm curl}_y\Big(\overline{Q^{-1}}\big({x\over\ep}\big)\overline U\big(t,x,{t\over \ep},{x\over\ep}\big)\Big)\Big)dxdt\\ \ecart\dis
=\int_0^T\hskip-3pt\int_{\RR^3}Q^{-1}\big({x\over\ep}\big)f_\ep\cdot \overline{U}\big(t,x,{t\over \ep},{x\over\ep}\big).\ea\eeq
Analogously, taking $\overline{V}(t,x,t/ \ep,x/\ep),$
as test function in the second equation in (\ref{maxwco}), we have
\beq\label{e2dtp}\ba{l}\dis -\int_{\RR^3}P\big({x\over\ep}\big)E_\ep^0\cdot \overline V\big(0,x,0,{x\over\ep}\big)dx\\ \ecart\dis-\int_0^T\hskip-3pt\int_{\RR^3}P\big({x\over\ep}\big)E_\ep\cdot \Big( \partial_t\overline V\big(t,x,{t\over \ep},{x\over\ep}\big)+{1\over\ep}\partial_s\overline V\big(t,x,{t\over \ep},{x\over\ep}\big)\Big)dxdt\\ \ecart\dis-\int_0^T\hskip-3pt\int_{\RR^3}Q^{-1}\big({x\over\ep}\big)B_\ep\cdot \Big( {\rm curl}_x\overline V\big(t,x,{t\over \ep},{x\over\ep}\big)+{1\over\ep}{\rm curl}_y\overline V\big(t,x,{t\over \ep},{x\over\ep}\big)\Big)dxdt\\ \ecart\dis
=\int_0^T\hskip-3pt\int_{\RR^3}g_\ep\cdot \overline{V}\big(t,x,{t\over \ep},{x\over\ep}\big).\ea\eeq\par
Taking into account that $(U,V)(t,x,.,.)$ belongs to ${\cal M}$ for a.e. $(t,x)\in (0,T)\times \RR^3$, we have
$$\ba{l}\dis \int_0^T\hskip-3pt\int_{\RR^3}Q^{-1}\big({x\over\ep}\big)B_\ep\cdot \Big(\partial_s\overline U\big(t,x,{t\over \ep},{x\over\ep}\big)+ {\rm curl}_y\overline V\big(t,x,{t\over \ep},{x\over\ep}\big)\Big)dxdt=0,
\ea$$
$$\int_0^T\hskip-3pt\int_{\RR^3}E_\ep\cdot \Big(\overline P\big({x\over\ep}\big)\partial_s\overline V\big(t,x,{t\over \ep},{x\over\ep}\big)- {\rm curl}_y\Big(\overline{Q^{-1}}\big({x\over\ep}\big)\overline U\big(t,x,{t\over \ep},{x\over\ep}\big)\Big)dxdt=0.$$
Therefore, adding (\ref{e1dtp}), (\ref{e2dtp}) we have that the terms which multiply $1/\ep$ cancel each other. By (\ref{conv2eda}), (\ref{conv2eso}), this allows us to pass to the limit in this sum to get
\beq\label{e3dtp}\ba{l}\dis -\int_{\RR^3}M_y\big(Q^{-1}\hat B^0\cdot \overline U(0,x,0,y)+P\hat E^0\cdot \overline V(0,x,0,y)\big)dx\\ \ecart\dis
-\int_0^T\hskip-3pt\int_{\RR^3}M_{sy}\big(Q^{-1}\hat B\cdot \partial_t\overline U+P\hat E\cdot\partial_t\overline V\big) dxdt\\ \ecart\dis
+\int_0^T\hskip-3pt\int_{\RR^3}M_{sy}\big(\hat E\cdot  {\rm curl}_x\big(\overline{Q^{-1}U}\big)-Q^{-1}\hat B\cdot {\rm curl}_x\overline V\big)dxdt\\ \ecart\dis
=\int_0^T\hskip-3pt\int_{\RR^3}M_{sy}\big(Q^{-1}\hat f\cdot\overline U+\hat g\cdot\overline V\big)dxdt.
\ea\eeq
By linearity and density, this equality holds replacing  $(U,V)$ by a couple $(U\varphi,V\varphi)$ with
$(U,V)$ an arbitrary element of ${\cal M}$ and $\varphi\in C^\infty([0,T];C^\infty_c(\RR^N))$, with $\varphi(T,x)=0$. This proves (\ref{pbli}).\par
The uniqueness of (\ref{pbli}) follows from Remark \ref{Rmcae} below (see also \cite{CCMM2}, Theorem B.1). As a consequence  we deduce that there is no need to extract a subsequence in (\ref{conv2eso}). \cqfd
\begin{Rem}\label{Recoda} Formally,  Theorem \ref{thcatwc} can be obtained 
  looking for an ansatz for the solution $(B_\ep,E_\ep)$ of (\ref{maxwco}) of the form
 \beq\label{anBe}B_\ep(t,x)\sim \hat B\Big(t,x,{t\over\ep},{x\over\ep}\Big)+\ep\hat B^1\Big(t,x,{t\over\ep},{x\over\ep}\Big)+\cdots\eeq
  \beq\label{anEe}E_\ep(t,x)\sim \hat E\Big(t,x,{t\over\ep},{x\over\ep}\Big)+\ep\hat E^1\Big(t,x,{t\over\ep},{x\over\ep}\Big)+\cdots\eeq
with $\hat B$, $\hat E$ almost-periodic in the microscopic variables $s=t/\ep$ $y=x/\ep$. To simplify, assume 
$$f_\ep(t,x)=\hat f\Big(t,x,{t\over\ep},{x\over\ep}\Big),\quad g_\ep(t,x)=\hat g\Big(t,x,{t\over\ep},{x\over\ep}\Big),\quad B^0_\ep(x)=\hat B^0\Big(x,{x\over\ep}\Big),\quad E^0_\ep(x)=\hat E^0\Big(x,{x\over\ep}\Big),$$
with 
\beq\label{conddiR}{\rm div}_y\hat f= {\rm div}_y\hat g=0,\quad {\rm div}_y\hat B^0= {\rm div}_y\hat E^0=0.\eeq
Replacing the ansatz for $B_\ep$, $E_\ep$ in the first and second equations in (\ref{maxwco}), 
and identifying the terms of equal power in $\ep$ we get
 \beq\label{prisif}\left\{\ba{l}\dis\partial_s\hat B+{\rm curl}_y\hat E=0\\ \ecart\dis P\partial_s\hat E+{\rm curl}_y(Q^{-1}\hat B)=0,\ea\right.\eeq
 \beq\label{segsif}\left\{\ba{l}\dis \partial_t\hat B+{\rm curl}_x\hat E+\partial_s\hat B^1+{\rm curl}_y\hat E^1=\hat f\\ \ecart\dis P\partial_t\hat E+{\rm curl}_x(Q^{-1}\hat B)+P\partial_s\hat E^1+{\rm curl}_y(Q^{-1}\hat B^1)=\hat g,\ea\right.\eeq
 \beq\label{inicohBhE} \hat B(0,x,0,y)=\hat B(x,y),\quad \hat E(0,x,0,y)=\hat E(x,y).\eeq
Taking the divergence in the first equation in (\ref{prisif}), we have that $\partial_s{\rm div}_y\hat B=0$ and therefore ${\rm div}_y\hat B={\rm div}_y M_s(\hat B)$. On the other hand, taking the mean value in $s$, and then the divergence in $y$, in the first equation in (\ref{segsif}), we get
$$\partial_t{\rm div}_yM_s(\hat B)+{\rm div}_y{\rm curl}_xM_s(\hat E)=0,$$
where using (\ref{prisif}), we have
$${\rm div}_y{\rm curl}_xM_s(\hat E)=-{\rm div}_x{\rm curl}_yM_s(\hat E)={\rm div}_x{\rm curl}_yM_s(\partial_s \hat B)=0.$$
Therefore, $\partial_t{\rm div}_y\hat B=0$, which combined with $\partial_s{\rm div}_y\hat B=0$ and
the first condition in (\ref{inicohBhE})  proves that ${\rm div}_y\hat B=0$. A similar reasoning also proves that ${\rm div}_y(P\hat E)=0$. Taking into account (\ref{prisif}) we conclude that as function of the variables $(s,y)$ the pair $(\hat B,\hat E)$ belongs to ${\cal M}$. Since this space has infinite dimension, in order to characterize $(\hat B,\hat E)$ we  return to (\ref{segsif}). To eliminate $\hat E^1$, $\hat B^1$ in this system we use 
$(Q^{-1}\hat U,\hat V)$ as test function, with $(U,V)$ valued in ${\cal M}$. Using that
$$M_{sy}\big((\partial_s\hat B^1+{\rm curl}_y\hat E^1)\cdot Q^{-1}\overline U+(P\partial_s\hat E^1+{\rm curl}_y(Q^{-1}\hat B^1))\cdot \overline V\big)=0,$$
we get (\ref{pbli}).
 \end{Rem}
 \begin{Rem} \label{Rmcae} To effectively compute the solution $(\hat B,\hat E)$ of (\ref{pbli}) we can use that by  Proposition \ref{ThcAu}, the functions $\hat B$, $\hat E$ admit the decomposition
\beq\label{deBE}\hat B=\sum_{\xi\in Y^3}\sum_{j=1}^\infty e^{i(\lambda^\xi_j s+\xi\cdot y)}\sum_{l=1}^{m_j^\xi} a^\xi_{jl}(t,x)\Phi^\xi_{jl}(y),\quad
\hat E=\sum_{\xi\in Y^3}\sum_{j=1}^\infty e^{i(\lambda^\xi_j s+\xi\cdot y)}\sum_{l=1}^{m_j^\xi} a^\xi_{jl}(t,x)\Psi^\xi_{jl}(y).\eeq
Taking in (\ref{pbli}) 
$$(U,V)=e^{i(\lambda^\xi_js+\xi\cdot y)}\big(\Phi^\xi_{jl},\Psi^\xi_{jl}\big),$$ 
we deduce that for every $\xi\in Y^3$ and every $j\geq 1$, the coefficients $a^\xi_{jl}$ are the solutions of the finite-dimensional hyperbolic linear system
\beq\label{sisthau}\left\{\ba{l}\dis\partial_t a^\xi_{jl}+\sum_{q=1}^{m^\xi_j}b^\xi_{jlq}\cdot\nabla_x a^\xi_{jq}=f_{jl}^\xi\\ \ecart\dis
(a^\xi_{jl})_{|t=0}= M_y\big(e^{-i\xi\cdot y}\big(Q^{-1}\hat B^0\cdot\overline {\Phi^\xi_{jl}}+P\hat E^0\cdot\overline {\Psi^\xi_{jl}}\big)\big),\ea\right.\qquad 1\leq l\leq m_j,\eeq
with
\beq\label{defjlx}f_{jl}^\xi=M_{sy}\big(e^{-i(\lambda^\xi_j s+\xi\cdot y)}\big(Q^{-1}\hat f\cdot\overline {\Phi^\xi_{jl}}+\hat g\cdot\overline {\Psi^\xi_{jl}}\big)\big),\eeq
\beq\label{defbjlq} b^\xi_{jlq}=\int_Y\big(\Psi_{jq}^\xi\cdot \overline{Q^{-1}\Phi^\xi_{jl}}-Q^{-1}\Phi^\xi_{jq}\cdot \overline{\Psi_{jl}^\xi}\big)dy.\eeq\par
Since by (\ref{condivd}), $(\hat B^0,\hat E^0)$ satisfies ${\rm div}_y\hat B^0={\rm div}_y(P\hat E^0),$ we can apply Proposition \ref{ThcAu} to deduce that $(\hat B^0,\hat E^0)$ can be decomposed as
$$(\hat B^0,\hat E^0)=\sum_{\xi\in Y^3}\sum_{j=1}^\infty  e^{i\xi\cdot y}\sum_{l=1}^{m_j^\xi}\alpha_{jl}^\xi\big(\Phi^\xi_{jl}(y), \Psi^\xi_{jl}(y)\big),$$
with 
$$\alpha^\xi_{jl}=M_{y}\Big(e^{-i\xi\cdot y}\big(Q^{-1}\hat B^0\cdot\Phi^\xi_{jl}+P\hat E^0\cdot \Psi^\xi_{jl}\big)\Big)=(a^\xi_{jl})_{|t=0}
.$$
This proves that
$$\big(\hat B(0,x,0,y),\hat E(0,x,0,y)\big)=\sum_{\xi\in Y^3}\sum_{j=1}^\infty \sum_{l=1}^{m_j^\xi} a^\xi_{jl}(0,x)e^{i\xi\cdot y}\big(\Phi^\xi_{jl}(y).\Psi^\xi_{ij}(y)\big)=\big(\hat B^0,\hat E^0\big),$$
and then it allows us to write the initial conditions for $\hat B$ and $\hat E$ in the simple form
\beq\label{coniBE} (\hat B)_{t=0,s=0}=\hat B^0,\quad (\hat E)_{t=0,s=0}=\hat E^0.\eeq
\end{Rem}
As a consequence of Theorem \ref{thcatwc}, we can now obtain the limit problem of (\ref{maxwco}). This is given by Theorem  \ref{thcod}. Before, we need 
\begin{Def}\label{DeQhPh}
 We define the homogenized matrices corresponding to $P$ and $Q$  by
\beq\label{dePQh} P_H\zeta=\int_YP(\zeta+\nabla Z_\zeta)\,dy,\quad Q_H\eta=\int_YQ(\eta+\nabla W_\eta)\,dy,\quad \quad \forall\,\zeta,\eta\in\RR^3,\eeq
with $W_\eta$, $Z_\zeta$ defined by (\ref{MUMV1}).\par
\end{Def}
\begin{Thm}\label{thcod} In the conditions of Theorem \ref{thcatwc}, defining
$f,g,B^0,(PE)^0$ by (\ref{convdd}) and $P_H$, $Q_H$ by Definition \ref{DeQhPh}, we have
\beq\label{convdBE} B_\ep\stackrel{\ast}\rightharpoonup B,\ \ E_\ep\stackrel{\ast}\rightharpoonup E\quad\hbox{ in }L^\infty(0,T;L^2(\RR^3))^3,\eeq
with $B,E$ the solutions of
\beq\label{Eclim}\left\{\ba{l}\dis B,E\in L^\infty(0,T;L^2(\RR^3))^3\\ \ecart\dis
\partial_tB+{\rm curl}\,E=f\ \hbox{ in }(0,T)\times \RR^N\\ \ecart\dis
P_H\partial_tE-{\rm curl}\, (Q^{-1}_HB)=g\ \hbox{ in }(0,T)\times \RR^N\\ \ecart\dis
B_{|t=0}=B^0,\ P_HE_{|t=0}=(PE)^0\ \hbox{ in }\RR^3.
\ea\right.\eeq
Moreover, the functions the functions $\hat B$, $\hat E$ in Theorem \ref{thcatwc} are related with $B$, $E$ by
\beq\label{relahBB} B=M_{sy}\hat B,\quad E=M_{sy}\hat E,\quad Q_H^{-1}B+\nabla_yW_{Q_H^{-1}B}=Q^{-1}M_s\hat B,\quad E+\nabla_yZ_E=M_sE,\eeq
with $W_{Q_H^{-1}B}$, and $Z_E$ are defined by (\ref{MUMV1}).
\end{Thm}\par\medskip\noindent
{\bf Proof.} Due to the relationship between weak convergence and two-scale convergence (see Theorem \ref{Thco2e}) and $B_\ep,E_\ep$ bounded in $L^\infty(0,T;L^2(\RR^3))^3$, we have that (\ref{convdBE}) holds with 
$B=M_{sy}\hat B$, $E=M_{sy}\hat E$. On the other hand, statement 5 in Proposition \ref{ThcAu} implies
\beq\label{camsBE} M_s\hat B= Q(\eta_B+\nabla_yW_{\eta_B}),\quad M_s\hat E=E+\nabla_yZ_E,\eeq
with $\eta_B\in L^\infty(0,T;L^2(\RR^3))^3$ and $W_{\eta_B}$, $Z_E$ defined by (\ref{MUMV1}). Namely, defintion (\ref{dePQh}) of $Q_H$ provides
\beq\label{cop0} B=\int_{Y^3}Q(\eta_B+\nabla_yW_{\eta_B})dy=Q_H\eta_B\Longrightarrow \eta_B=Q_H^{-1}B.\eeq
Taking in (\ref{pbli}) $U=Q(\eta+\nabla_yW_\eta)$, $V=\zeta+\nabla_yZ_\zeta$, with $\eta,\zeta\in \RR^3$, $W_\eta,Z_\zeta\in H^1_\sharp(Y^3)$, related by (\ref{MUMV1}), we get
\beq\label{cop1}\ba{l}\dis \partial_t\int_{Y^3}\Big(Q(\eta_B+\nabla_yW_{\eta_B})\cdot \overline{(\eta+\nabla_yW_\eta)}+P(E+\nabla_yZ_E)\cdot \overline{(\zeta+\nabla_yZ_\zeta)}\Big)dy\\ \ecart\dis
+{\rm div}_x\int_{Y^3}\Big((E+\nabla_yZ_E)\times \overline{(\eta+\nabla_yW_\eta)}-(\eta_B+\nabla_yW_{\eta_B})\times \overline{(\zeta+\nabla_yZ_{\zeta})}\Big)dy\\ \ecart\dis
=\int_{Y^3}\Big(M_s\hat f\cdot \overline{(\eta+\nabla_yW_\eta)}+M_s\hat g\cdot \overline{(\zeta+\nabla_yZ_{\zeta})}\Big)dy,
\ea\eeq
combined with the initial conditions
\beq\label{cop2}\ba{l}\dis \int_{Y^3}\Big(Q(\eta_B+\nabla_yW_{\eta_B})\cdot \overline{(\eta+\nabla_yW_\eta)}+P(E+\nabla_yZ_E)\cdot \overline{(E+\nabla_yZ_\zeta)}\Big)dy\\ \ecart\dis =\int_{Y^3}\Big(\hat B^0\cdot \overline{(\eta+\nabla_yW_\eta)}+P\hat E^0\cdot\overline{(\zeta+\nabla_yZ_{\zeta})}\Big)dy,\quad \hbox{ for }t=0.\ea
\eeq
\par
The equations (\ref{MUMV1}) satisfied by $W_{\eta_B}$ and $Z_E$, definition (\ref{dePQh}) of $P_H$, $Q_H$ and (\ref{cop0}) give
$$\int_{Y^3}\Big(Q(\eta_B+\nabla_yW_{\eta_B})\cdot \overline{(\eta+\nabla_yW_\eta)}+P(E+\nabla_yZ_E)\cdot \overline{(\zeta+\nabla_yZ_\zeta)}\Big)dy=B\cdot \overline\eta+P_HE\cdot \overline\zeta.$$
Using also that by  (\ref{cop0})
$$\int_{Y^3}\Big((E+\nabla_yZ_E)\times \overline{(\eta+\nabla_yW_\eta)}-(\eta_B+\nabla_yW_{\eta_B})\times \overline{(\zeta+\nabla_yZ_{\zeta})}\Big)dy=E\times\overline \eta-Q_H^{-1}B\times\overline\zeta,$$
that ${\rm div}_y\hat f={\rm div}_y\hat g=0$, and that definitions (\ref{limded})  of $f$ and $g$ imply
$$\int_{Y^3}\Big(M_s\hat f\cdot \overline{(\eta+\nabla_yW_\eta)}+M_s\hat g\cdot \overline{(\zeta+\nabla_yZ_{\zeta})}\Big)dy=f\cdot\overline\eta+g\cdot\overline\zeta,$$
we get that (\ref{cop1}) reduces to
\beq\label{cop3}\partial_t\big(B\cdot \overline\eta+P_HE\cdot \overline\zeta\big)+{\rm div}_x\big(E\times\overline \eta-Q_H^{-1}B\times\overline\zeta\big)=f\cdot\overline\eta+g\cdot\overline\zeta\ \hbox{ in }(0,T)\times\RR^3.\eeq\par
Using the equations satisfied by $W_{\eta_B}$ and $Z_E$, that ${\rm div}_y\hat B^0={\rm div}_y(P\hat E^0)=0$, and $B^0=M_y\hat B^0$, $(PE)^0=M_y(P\hat E^0)$, we also have that (\ref{cop2}) can be written as
\beq\label{cop4}B\cdot\overline{\eta}+P_HE\cdot\overline\zeta=\hat B^0\cdot\overline{\eta}+M_y(P\hat E^0)\cdot\overline\zeta\ \hbox{ in }t=0.\eeq
Equalities (\ref{cop3}) and (\ref{cop4}) with $\eta$ and $\zeta\in\RR^3$ arbitrary, are equivalent to (\ref{Eclim}).\par

Statement (\ref{relahBB}) follows from $B=M_{sy}\hat B$, $E=M_{sy}\hat E$, and (\ref{camsBE}) with $\eta_B=Q_H^{-1}B$.\cqfd
\begin{Rem} \label{cobp}
Theorem \ref{thcod} is known (see e.g. \cite{SaPa}, Chapter 7, \cite{Wei}). It  can be easily obtained using for example the two-scale convergence theory for periodic functions (\cite{All}, \cite{Ngu}). Assuming that the right-hand sides $f_\ep$, $g_\ep$ are of the form 
\beq\label{rsndt} f_\ep(t,x)=\hat f(t,x,x/\ep),\quad g_\ep(t,x)=\hat g(t,x,x/\ep), \eeq
with $f,g$ periodic (or almost periodic) in the quick variable $y=x/ep$, it can be formally obtained seeking an ansatz for $(B_\ep, E_\ep)$ of the form 
\beq\label{anBeb}B_\ep(t,x)\sim \check B\Big(t,x,{x\over\ep}\Big)+\ep \check B^1\Big(t,x,{x\over\ep}\Big)+\cdots,\quad E_\ep(t,x)\sim  \check E\Big(t,x,{x\over\ep}\Big)+\ep \check E^1\Big(t,x,{x\over\ep}\Big)+\cdots.\eeq
Observe that contrarily to the ansatz in Remark \ref{Recoda}, we are assuming that the functions $B_\ep$ and $E_\ep$ do not oscillate in the time variable. This may seem natural because if $f_\ep$ and $g_\ep$ satisfy (\ref{rsndt}), then neither coefficients nor the right-hand sides in (\ref{maxwcoU}) oscillate in the time variable. The corresponding functions $\check B$ and $\check E$ must be obtained from $\hat E$ and $\hat B$ by $\check B=M_s\hat B$ and $\check E=M_s\hat E$. Thus, using (\ref{relahBB}) we get
\beq\label{corcl}\check B=Q\big(Q_H^{-1}B+\nabla_yW_{Q_H^{-1}B}\big),\quad \check E=E+\nabla_yZ_E.\eeq
As it is usual in classical homogenization these functions depend locally on the limits $B$ and $E$ of the sequences $B_\ep$ and $E_\ep$ in the sense that for a.e. $(t,x)\in (0,T)\times\RR^3$, the values of $\check B(t,x,.)$ and  $\check E(t,x,.)$ only depend on $B(t,x)$ and $E(t,x)$. If $\check B$ and $\check E$ effectively provide an approximation in the strong topology of $L^2(0,T;L^2(\RR^3))^3$ for $B_\ep$ and $E_\ep$ (corrector result) then we must have $\hat B=\check B$ and $\hat E=\check E$ which is equivalent to have in (\ref{deBE}) that all the functions $a_{jl}^\xi$ vanish except those corresponding to $\xi=0$ and $\lambda_j^\xi=0$. Taking into account (\ref{sisthau}) where, thanks to (\ref{rsndt}), the coefficients $f^\xi_{jl}$ vanish, we deduce that  $(\hat B^0,\hat E^0)$ must be orthogonal for the scalar product given by (\ref{proeL2}) to all the functions $(\Phi^\xi_{jl},\Psi^\xi_{jl})$ corresponding to $\xi\not=0$ or $\xi=0$, $\lambda^\xi_j\not =0$, i.e. by statement 5 in Proposition \ref{ThcAu} they must satisfy
\beq\label{conico}\hat B^0=Q\big(Q_H^{-1}B^0+\nabla_yW_{Q_H^{-1}B^0}\big),\quad \hat E^0=E^0+\nabla_yZ_{E^0},\eeq
This holds if we assume  that $B_\ep^0$ and $E_\ep^0$ satisfy (in addition to (\ref{consd})) that
\beq\label{dabp}{\rm curl}\big(Q\big({x\over \ep}\big)^{-1}B^0_\ep\big)\ \hbox{ is compact in }H^{-1}(\RR^3)^3,\quad {\rm div}\big(P\big({x\over \ep}\big)E^0_\ep\big)\ \hbox{ is compact in }H^{-1}(\RR^3),\eeq
which would be the equivalent in our context to the well posed conditions in \cite{FrMu} for the homogenization of the wave equation. However (\ref{conico}) does not hold even if we assume that $B^0_\ep$, $E^0_\ep$ do not depend on $\ep$. This means that contrarily to assumption (\ref{anBeb}) the solutions of (\ref{maxwco}) oscillate in time even if the right-hand sides do not. \par
Observe that while $\check B$, $\check E$ can be locally obtained from $E$ and $B$, the remaining terms in decomposition (\ref{deBE}) of $\hat E$, $\hat B$, which solve (\ref{sisthau}), depend non-locally on the initial conditions and the right-hand sides.
\end{Rem}
The following theorem combined with Proposition \ref{Procf} provides now a corrector result for the solution of (\ref{maxwco}).
\begin{Thm} \label{Thcorr} In the conditions of Theorem \ref{thcatwc}. If the  convergences in (\ref{conv2eda}) hold in the strong two-scale convergence sense, then convergences (\ref{conv2eso}) also hold in the strong two-scale convergence sense.
\end{Thm}
\par\noindent
{\bf Proof.} Integrating twice in time (\ref{EnId}) applied to (\ref{maxwco}), we get
$$\ba{l}\dis {1\over 2}\int_0^T\hskip-3pt\int_{\RR^3}\Big(Q^{-1}\big({x\over\ep}\big)B_\ep\cdot\overline{B_\ep}+P\big({x\over\ep}\big)E_\ep\cdot\overline{E_\ep}\Big)dxdt\\ \ecart\dis
={T\over 2}\int_{\RR^3}\Big(Q^{-1}\big({x\over\ep}\big)B^0_\ep\cdot\overline{B^0_\ep}+P\big({x\over\ep}\big)E^0_\ep\cdot\overline{E^0_\ep}\Big)dx\\ \ecart\dis+
Re\int_0^T\hskip-3pt\int_{\RR^3}(T-t)\Big(Q^{-1}\big({x\over\ep}\big)f_\ep\cdot\overline{B_\ep}+g_\ep\cdot\overline{E_\ep}\Big)dxdt.\ea$$
Since we are assuming that the convergences in (\ref{conv2eda}) hold in the strong two-scale sense, we can use (\ref{conv2eso}) to deduce
\beq\label{e1thco}\ba{l}\dis \lim_{\ep \to 0}{1\over 2}\int_0^T\hskip-3pt\int_{\RR^3}\Big(Q^{-1}\big({x\over\ep}\big)B_\ep\cdot\overline{B_\ep}+P\big({x\over\ep}\big)E_\ep\cdot\overline{E_\ep}\Big)dxdt\\ \ecart\dis={T\over 2}\int_{\RR^3}M_y\big(Q^{-1}\hat B^0\cdot\overline{\hat B^0}+P\hat E^0\cdot\overline{\hat E^0}\big)dx\\ \ecart\dis
+
Re\int_0^T\hskip-3pt\int_{\RR^3}(T-t)M_{sy}\Big(Q^{-1}\hat f\cdot\overline{\hat B}+\hat g\cdot\overline{\hat E}\Big)dxdt.\ea\eeq
On the other hand, taking in (\ref{pbli}) $U=\hat B$ and $V=\hat E$, and taking the real part,  we get (the functions $\hat B$ and $\hat E$ are not the necessary regular but the result can be obtained by using an approximation of them) 
$$
 {1\over 2}{d\over dt}\int_{\RR^3} M_{sy}\big(Q^{-1}\hat B\cdot \overline {\hat B}+P\hat E\cdot\overline {\hat E}\big) dxdt=Re\int_0^t\hskip-3pt\int_{\RR^3}M_{sy}\big(Q^{-1}\hat f\cdot\overline {\hat B}+\hat g\cdot\overline {\hat E}\big)dxdt.$$
Integrating twice this inequality and using the boundary conditions in (\ref{pbli}) combined with (\ref{coniBE}), we also get
$$\ba{l}\dis {1\over 2}\int_0^T\hskip-3pt\int_{\RR^3}M_{sy}\big(Q^{-1}\hat B\cdot \overline {\hat B}+P\hat E\cdot\overline {\hat E}\big) dxdt\\ \ecart\dis={T\over 2}\int_{\RR^3}M_s\big(Q^{-1}\hat B^0\cdot\overline{\hat B^0}+P\hat E^0\cdot\overline{\hat E^0}\Big)dx+
Re\int_0^T\hskip-3pt\int_{\RR^3}(T-t)M_{sy}\Big(Q^{-1}\hat f\cdot\overline{\hat B}+\hat g\cdot\overline{\hat E}\Big)dxdt.\ea$$
By (\ref{e1thco}), this proves the convergence of the energies
$$ \lim_{\ep \to 0}{1\over 2}\int_0^T\hskip-3pt\int_{\RR^3}\Big(Q^{-1}\big({x\over\ep}\big)B_\ep\cdot\overline{B_\ep}+P\big({x\over\ep}\big)E_\ep\cdot\overline{E_\ep}\Big)dxdt={1\over 2}\int_0^T\hskip-3pt\int_{\RR^3}M_{sy}\big(Q^{-1}\hat B\cdot \overline {\hat B}+P\hat E\cdot\overline {\hat E}\big) dxdt,$$
and then the strong two-scale convergence of $B_\ep$ and $E_\ep$. \cqfd

\section*{Aknowledgments} The work of J. Casado-D\'{\i}az has been partially supported by the project PID2020-116809GB-I00
of the {\it Ministerio de Ciencia e Innovaci\'on}  of the goverment of Spain. The work of N. Khedhiri, is supported by the Tunisian ministry of higher education. The work of M.-L. Tayeb is supported by Tamkeen under the NYU Abu Dhabi Research Institute grant of the
center SITE.

\par\vfill\eject


\begin{thebibliography}{20}

\bibitem{All} {\sc G. Allaire.} {\it Homogenization and two-scale convergence}. SIAM J. Math. Analysis 23 (1992), 1482-1518.

\bibitem{AlFr} {\sc G. Allaire, L. Friz.} {\it Localization of high-frequence waves propagating in a locally periodic medium}. Proc. Roy. Soc. Edinburgh 140A (2010), 897-926.

\bibitem{APR} {\sc G. Allaire, M. Palombaro, J. Rauch}. {\it Diffractive behavior of the wave equation in periodic media: weak convergence analysis}.
Annali di Matematica  188 (2009), 561-589.

\bibitem{BFM} {\sc S. Brahim-Otsmane, G.A. Francfort, F. Murat}. {\it Correctors for the homogenization of the wave and heat equations}. J. Math. Pures Appl. 71 (1992), 197-231.

\bibitem{BrLe} {\sc M. Brassart, M. Lenczner}. {\it A two scale model for the periodic homogenization of the wave equation}. J. Math. Pur. Appl. 93 (2010), 474-517.

\bibitem{BrCa} {\sc M. Briane, J. Casado-D\'{\i}az.} {\it Increase of mass and nonlocal effects in the homogenization of magneto-elastodynamics problems.} Calc. Var.  and PDE, 60, 5 (2021), paper 163, 39 pp.

\bibitem{BLP} {\sc A. Bensoussan, J.L. Lions, G. Papanicolau.} {\it Asymptotic analysis for periodic structures.} North-Holland, Amsterdam-New York-Oxford, 1978.

\bibitem{CCMM1} {\sc J. Casado-D\'{\i}az, J. Couce-Calvo,  F. Maestre, J.D. Mart\'{\i}n-G\'omez}. {\it Homogenization and corrector for the wave equation with discontinuous coefficients in time.} J. Math. Anal. Appl. 379 (2011),  664--681.

\bibitem{CCMM2}{\sc J. Casado-D\'{\i}az, J. Couce-Calvo, F. Maestre, J.D. Mart\'{\i}n-G\'omez}. {\it Homogenization and correctors for the wave equation with periodic coefficients}. Math. Mod. Meth. Appl. Sci. 24 (2014), 1343-1388.

\bibitem{CCMM3}{\sc J. Casado-D\'{\i}az, J. Couce-Calvo, F. Maestre, J.D. Mart\'{\i}n-G\'omez}. {\it A corrector for a wave problem with periodic coefficients in a 1D bounded domain}. ESAIM: Cont. Optim. Calc. Var. 21 (2015), 465-486.

\bibitem{CaGa} {\sc J. Casado-D\'{\i}az, I. Gayte}. {\it A general compactness result and its application the two-scale convergence of almost periodic functions}. C. R. Acad. Sci. Paris I 323 (1996), 329-334.

\bibitem{CaGa2}{\sc J. Casado-D\'{\i}az, I. Gayte}. {\it The two-scale convergence method applied to generalized Besicovitch spaces}. Proc. Roy. Soc. London A 458 (2002), 2925-2946.

\bibitem{CaMa} {\sc J. Casado-D\'{\i}az, F. Maestre}. {\it An elliptic equation in an unbounded cylinder. Applications to the behavior of a wave in a thin beam with boundary conditions.} Rev. Mat. Complutense 32 (2019), 681-730.

\bibitem{CoSp} {\sc F. Colombini, S. Spagnolo}. {\it On the convergence of solutions of hyperbolic equations.} Commun. Partial Differential Equations 3 (1978), 77--103.

\bibitem{FrMu}  {\sc G.A. Francfort, F. Murat}. {\it Oscillations and energy densities in the wave equation}. Comm. Partial Differential Equations
17 (1992), 1785-1865.

\bibitem{HMC} {\sc D.  Harutyunyan, G.W. Milton, R.V. Craster}. {\it High-frequency homogenization for travelling waves in periodic media.} Proc. A 472 (2016), 18 pp.

\bibitem{Mur}  {\sc F. Murat, L. Tartar.} {\it $H$-convergence}. In: {\sc L. Cherkaev, R.V. Kohn} (Eds.). {\it Topics in the Mathematical Modelling of Composite Materials}. Progress in Nonlinear Differential Equations and Their Applications, vol. 31, Birk\"auser, Boston, 1998, pp. 21--43.

\bibitem{Ngu} {\sc G. Nguetseng}. {\it A general convergence result for a functional related to the theory of homogenization}. SIAM J. Math. Analysis 20 (1989), 608-623.
\bibitem{Ngu2} {\sc G. Nguetseng.} {\it Homogenization structures and applications I}. Z. Anal. Anwen 22 (2003) 73--107.

\bibitem{NgLeBr} {\sc T.T. Nguyen, M. Lenczner, M. Brassart.} {\it Homogenization of the one-dimensional wave equation.} In {\sc A. Abdulle, S. Deparis, D. Kressner, F. Nobile, M. Picasso} (Eds.). {\it Numerical Mathematics and advances applications - ENUMATH 2013}. Lect. Notes Comput. Sci. Eng., 103, Springer, Cham, 2015, pp. 377--385.

\bibitem{Spa} {\sc S. Spagnolo.} {\it Sulla convergenza di soluzioni di equazioni paraboliche ed ellittiche.} Ann. Scuola Norm. Sup. Pisa Cl. Sci. 22 (1968), 571--597.

\bibitem{SaPa} {\sc E. S\'anchez-Palencia.} {\it Non-Homogeneous Media and Vibration Theory.} Lect. Notes in Physics 127. Springer-Verlag, Berlin-Heidelberg-New-York, 1980.

\bibitem{Tar} {\sc L. Tartar.} {\it The General Theory of Homogenization: A Personalized Introduction.} Lect. Notes of the Unione Matematica Italiana, Springer-Verlag, Berlin, Heidelberg, 2009.

\bibitem{Wei} {\sc N. Wellander.} {\it Homogenization of the Maxwell equations: Case I. Linear theory.} Appl. Math. 46 (2001), 29--51.
\end{thebibliography}
\end{document}